\documentstyle[12pt]{article}
\begin{document}
\title{ Poisson measures for topological groups
and their representations.}
\author{ S.V. Ludkovsky}
\date{ 01 December 1998
\thanks{ Mathematics subject classification (1991 Revision) 
22A10, 43A05 and 46B.} }
\maketitle
\section{Introduction}
\par In articles \cite{luumn96,lutmf,luumnad,luihlgr,luihlgna,luihgdsw,lupr202,
lupr215,lupr218,shav}
Gaussian quasi-invariant measures on groups of diffeomorphisms and loop
groups $G$ relative to dense subgroups $G'$ were constructed.
In the non-Archimedean case the wider class of measures was investigated,
than in the real case.
The cases of Riemann and non-Archimedean manifolds 
were considered. There are few approaches for the construction of irreducible
unitary representations. In articles 
\cite{kos,lutmf,luihlgr,luihlgna,luihgdsw} representations
of dense subgroups $G'$
associated with quasi-invariant measures on the entire groups
were considered. In articles \cite{hira,luseamb,lurimut}
irreducible representations of groups of diffeomorphisms 
$Diff(M)$ associated with measures on specific 
subsets of the unital type of products $M^{\bf N}$ of the manifolds $M$
were investigated.
In the publications \cite{shim,vgg}
irreducible unitary representations of groups of diffeomorphisms
associated with real-valued Poisson measures
on products of real manifolds were studied.
\par This article is related with unitary representations of $G'$
associated with Poisson measures
on $G^{\bf N}$ and uses quasi-invariant measures on $G$ from the previous works.
Several groups are considered: (1) (a) diffeomorphisms and (b) loop 
groups of real
manifolds, (2) (a) diffeomorphisms and (b) loop groups of non-Archimedean 
manifolds over local fields.
Besides these four cases further the fifth and the sixth cases 
are considered:
for (3) (a) real and (b) non-Archimedean groups of diffeomorphisms $Diff(M)$
representations associated with
Poisson measures on configuration spaces $\Gamma _M$ contained
in products of manifolds $M^{\bf N}$ are investigated.
The case (3) (a) for real locally compact $M$ was considered in 
\cite{shim,vgg}. Here the cases of infinite-dimensional Banach
manifold $M$ (3) (a), non-Archimdean locally compact
and non-locally compact Banach manifolds (3) (b) are investigated.
For this quasi-invariant measures on $M$ relative to $Diff(M)$
from \cite{luseamb,lurimut} are considered.
Henceforth real-valued measures are considered.
In \S 2 necessary Poisson measures are considered, 
definitions and notations are
given. In \S 3 irreducible unitary representations are considered.
Certainly not all results from \cite{shim,vgg} can be transferred onto the
cases considered here, moreover, there were necessary strong changes
in many definitions, proofs and formulations of the theorems.
\par It is necessary to note that the theory of representations
of non-locally compact groups differ substantially from that
of locally compact groups. For example, irreducible unitary representations
of locally compact Abelian groups are one-dimensional, that is, characters.
But for non-locally compact Abelian groups there are infinite-dimensional 
irreducible unitary representations, which are even regular representations.
It was shown in \cite{banas,gelvil} that 
there are infinite-dimensional topological
vector spaces $E$ and dense nuclear additive subgroups $E'$ such that
$E'$ are linear subspaces and quasi-invariant measures $\mu $ on $E$ 
relative to $E'$ exist such that associated
with them regular representations in the Hilbert space $L^2(E,\mu ,{\bf C})$
are irreducible. The existence of such irreducible representations
is even despite of 
the fact that projections $\mu _J$ of $\mu $ on one-dimensional 
subspaces $J$ are equivalent with the Haar measures on $J$. This shows that
non-locally compact case is more complicated than it may be supposed at the
first glance.
Also for definite groups $G$ of diffeomorphisms and loops
of definite real and non-Archimedean manifolds there are quasi-invariant
measures $\mu $ on $G$ relative to dense subgroups $G'$ such that
associated with them regular unitary representations are irreducible
\cite{kos,lutmf,luihlgr,luihlgna,luihgdsw}.
Such difference is caused by the existence of $C^*$-algebras associated with
the Haar measures on locally compact groups \cite{hew}, 
but no any $C^*$-algebra can be directly associated with a non-zero
quasi-invariant measure on a non-locally
compact group relative to a dense subgroup $G'$.
Certainly, results on irreducibility of regular representations
of infinite-dimensional topological groups $G'$ depend stronlgy
on quasi-invariant measures $\mu $ on $G$ and a structure of $G$, 
where $G'$ is dense in $G$.
\section{Poisson measures.}
\par {\bf 2.1. Note.} Let $X$ denotes a manifold $M$ for a group of
diffeomorphisms $G=G(M)$ or the group $G$ itself, where 
$M$ is the $C^{\infty }$-manifold over $\bf R$ or an analytic
manifold over a local field
and $G$ is the loop group
or the diffeomorphisms group as in the cited in \S 1 papers.
Classes of smoothness of the groups $G$ and $G'$ 
are considered to be not less than $C^1$. The groups of diffeomorphisms
$G$ for the real $C^{\infty }$-manifold $M$ are denoted 
$Diff^t_{\beta ,\gamma }(M)$ with
$\infty \ge t\ge 1$, $\beta \ge 0$, $\gamma \ge 0$;
the loop groups $G$ for the real $C^{\infty }$-manifold $N$ 
are denoted $(L^mN)_{\gamma ,Y}$ with $m+5< \gamma \le \infty $,
also another classes of smoothness and non-Archimedean groups
and manifolds were considered (see theorem 3.4 \cite{lupr202}
and also \cite{ebi,lutmf,luihlgr,luihlgna,lurimut}).
It was proved earlier, that $G$ itself is the $C^{\infty }$-manifold
(in the case of the real group of diffeomorphisms for finite-dimensional
Riemann manifolds $M$ see also \cite{bao,ebi}).
Moreover, in the non-Archimedean case $M$ and $G$ have structures of 
the analytic manifolds with clopen disjoint charts.
Clearly, $G$ itself is not locally compact, since $G$ considered
as the manifold is infinite-dimensional over the corresponding field.
When $X=M$ let us suppose,
that $X$ is the Banach non-compact manifold. In the non-Archimedean 
case it has embedding into the Banach space $Z$ over the same local field
$\bf L$ due to the partition of $M$ into disjoint union of balls, so
an atlas of $M$ is supposed to be analytic and it has automatically
foliated structure \cite{luum985,lupr180}. In the real case
it is supposed that $M$ has a foliated structure with finite-dimensional
submanifolds $M_n\subset M_{n+1}$ for each $n\in \bf N$ and
$\bigcup_{n\in \bf N}M_n$ is dense in $M$, where $dim_{\bf R}M_n=k_n<\infty $
\cite{lurimut,lupr202,lupr218}.
\par We remind the definition of the configuration space
from \cite{shim} and also consider the ultrametric case of $X$.
This means that a metric $d$ in $X$ satisfies the ultrametric
inequality $d(x,y)\le \max (d(x,z), d(y,z))$ for each $x, y, z \in X$.
\par Let $K$ be a complete separable metric space with a metric $d$,
that is, $X$ is a Polish space.
In the ultrametric case this implies that its topological great inductive
dimension is zero: $Ind(K)=0$ \cite{eng}.
Let $d^n_K(x,y):=\sum_{i=1}^nd(x_i,y_i)$ in the real case
and $d(x,y):=\max_{1\le i \le n} d(x_i,y_i)$ in the non-Archimedean case
be a metric in $K^n$, where $x=(x_i: i=1,...,n)\in K^n$, $x_i\in K$.
Put ${\tilde K}^n:=(x\in K^n: x_i\ne x_j$
for each $i\ne j)$. Supply ${\tilde K}^n$ with a metric $\delta ^n_K(x,y):=
d^n_K(x,y)/[d^n_K(x,y)+d^n_K(x,({\tilde K}^n)^c)+d^n_K(y,({\tilde K}^n)^c)]$
in the real case and $\delta ^n_K(x,y):=d^n_K(x,y)/[\max (d^n_K(x,y),
d^n_K(x,({\tilde K}^n)^c), d(y,({\tilde K}^n)^c)]$
in the non-Archimedean case,
where $A^c:=K^n\setminus A$ for a subset $A\subset K^n$.
Then $({\tilde K}^n, \delta ^n_K)$ is the Polish space. Moreover,
if $(K,d)$ is ultrametric, then $({\tilde K}^n, \delta ^n_K)$
is ultrametric. Let also $B^n_K$ denotes the collection of all
$n$-point subsets of $K$. Then the metric $\delta ^n_K$ is equivalent with
the following metric $d^{(n)}_K(\gamma ,\gamma '):=\inf_{\sigma \in
\Sigma _n}d^n_K((x_1,...,x_n), (x'_{\sigma (1)},...,x'_{\sigma (n)}))$,
where $\Sigma _n$ is the symmetric group of $(1,...,n)$, $\sigma \in \Sigma _n$,
$\sigma : (1,...,n)\to (1,...,n)$; $\gamma , \gamma ' \in B^n_K$.
For each subset $A\subset K$ a number mapping $N_A: B^n_K\to \bf N_o$
is defined by the following formula: $N_A(\gamma ):=card(\gamma \cap A)$,
where ${\bf N}:=\{ 1,2,3,... \} $, ${\bf N_o}:=\{ 0,1,2,3,... \} $.
Evidently, $d^{(n)}_K$ is the ultrametric, if $d^n_K$ is the ultrametric.
It remains to show, that $\delta ^n_K$ is the ultrametric for the ultrametric
space $(K,d)$.  For this we mention, that $(i)$ $\delta ^n_K(x,y)>0$, when
$x\ne y$, and $\delta ^n_K(x,x)=0$. $(ii)$ $\delta ^n_K(x,y)=\delta ^n_K(y,x)$,
since this symmetry is true for $d^n_K$ and for $[*]$ in the denumerator
in the formula defining $\delta ^n_K$.  To prove $(iii)$ $\delta ^n_K(x,y)\le
\max (\delta ^n_K(x,z), \delta ^n_K(z,y))$ we consider the case
$\delta ^n_K(x,z)\ge \delta ^n_K(y,z)$, hence it is sufficient to show, that
$\delta ^n_K(x,y)\le \delta ^n_K(x,z)$.
Let $(a)$ $d^n_K(x,z)\ge \max (d^n_K(z,({\tilde K}^n)^c), d^n_K(x,({\tilde
K}^n)^c))$, then $\delta ^n_K(x,z)=1$, hence $\delta ^n_K(x,y)\le \delta
^n_K(x,z)$, since  $\delta ^n_K(x,y)\le 1$
for each $x, y\in {\tilde K}^n$. Let $(b)$ $d^n_K(x,({\tilde K}^n)^c)>
\max (d^n_K(x,z), d^n_K(z,({\tilde K}^n)^c))$, then $\delta ^n_K(x,z)=
d^n_K(x,z)/d^n_K(x,({\tilde K}^n)^c)\le 1$.
Since $d^n_K(z,A):=\inf_{a\in A}d^n_K(z,a)$, then $d^n_K(z,({\tilde K}^n)^c)
\le \max (d^n_K(y,({\tilde K}^n)^c), d^n_K(y,z))$. If $d^n_K(x,z)<
d^n_K(z,({\tilde K}^n)^c)$ and $d^n_K(x,y)\le d^n_K(x,z)$, then
$d^n_K(z,({\tilde K}^n)^c)\le d^n_K(x,({\tilde K}^n)^c)$. Hence
$d^n_K(x,y) \max (d^n_K(x,z), d^n_K(x, ({\tilde K}^n)^c), d^n_K(z,({\tilde
K}^n)^c))$ $\le d^n_K(x,z)\max (d^n_K(x,y), d^n_K(x,({\tilde K}^n)^c),
d^n_K(y,({\tilde K}^n)^c))$. With the help of $(ii)$ the remaining cases
may be lightly written. 
\par The Borel $\sigma $-field of $B^n_K$ is denoted by $Bf(B^n_K)$.
If $<S,{\sf L},m>$ is the measure space, then its completion relative to
$m$ is denoted $Af(S,m)$, where $S$ is a set, ${\sf L}$ is a $\sigma $-algebra
of subsets of $S$, $m$ is a real non-negative $\sigma $-finite
measure on ${\sf L}$.
That is, the $\sigma $-algebra $Af(S,m)$ contains all subsets $A\subset B$
of $B\in \sf L$ for which $m(B)=0$.
In the non-Archimedean case the valuation group ${\Gamma '}_{\bf L}
:=\{ |x|_{\bf L}:
0\ne x\in {\bf L} \} $ of the local field $\bf L$ is discrete in
$(0,\infty )$, hence subsets $U_{\epsilon }(y):=\{ x\in K: 
d(x,y)<\epsilon \} $ are clopen (closed and open simultaneously)
in $K=X$. Therefore, in the non-Archimedean case
lemmas 1.1 and 1.2 from \cite{shim} have the following stronger forms.
\par {\bf 2.2. Lemma.} {\it  For an ultrametric space $(X,d)$ from \S 2.1
if $U$ is a clopen set in $X$, then $\{ \gamma : N_U(\gamma )\ge l \} $
is also clopen in $X$ for each $l\in \bf N_o$.}
\par {\bf 2.3. Lemma.} {\it For an ultrametric space $(X,d)$ from \S 2.1
and each $\epsilon >0$ and each $\gamma \in B^n_X$
there exists a clopen subset $O_{\epsilon }(\gamma )$ which belongs to
the smallest $\sigma $-algebra $\sf B$ for which functions $N_B$ are measurable
such that $\gamma \in O_{\epsilon }(\gamma )\subset \{ \gamma ':
d^{(n)}_X(\gamma ,\gamma ')<\epsilon \} $.}
\par {\bf Proof.} For $\gamma =\{ x_1,...,x_n \} $ take $\eta \in 
{\Gamma '}_{\bf
L}$ such that $\epsilon >\eta  >0$ and $U_{\eta p^{-n}}(x_i)\cap U_{\eta
p^{-n}}(x_j)=\emptyset $ for each $i\ne j$. Put $O_{\epsilon }(\gamma )$ 
$:=
\{ \bigcap_{i=1}^n\{ \gamma ': card(\gamma '\cap U_{\eta p^{-n}}(x_i))
\ge 1 \} $, where $1<p\in \Gamma _{\bf L}$, $p^{-1}=|\pi _{\bf L}|_{\bf L}$,
$B({\bf L},0,1^{-})=\pi _{\bf L}B({\bf L},0,1)$, $B(Y,x,r):=\{
z\in Y: d_Y(x,z)\le r \} $, $B(Y,y,r^{-}):=\{ z\in Y: d_Y(y,z)<r \} $
for an ultrametric space $Y$ with an ultrametric $d_Y$.
\par {\bf 2.4. Notes and definitions.} Then theorems 1.1 and 1.2 
from \cite{shim} are also true for all cases
considered here. For this we mention, that as usually let 
$B_K:=\bigoplus_{n=0}^{\infty }B^n_K$,
where $B^0_K:=\{ \emptyset \} $ is a singleton. 
Since $X$ from \S 2.1 is not compact,
then there exists an increasing sequence of subsets $K_n\subset X$
such that $X=\bigcup_nK_n$ and $K_n$ are Polish spaces in the induced
topology from $X$. Moreover, $K_n$ can be chosen clopen in $X$
in the non-Archimedean case. Then the following space $\Gamma _X:=
\{ \gamma : \gamma \subset X$ and $card(\gamma \cap K_n)<\infty $
for each $n \} $ is called the configuration space and it is isomorphic
with the projective limit $pr-\lim \{ B_{K_n}, \pi ^n_m, {\bf N} \} $,
where $\pi ^n_m(\gamma _m)=\gamma _n$ for each $m>n$ and $\gamma _n
\in B_{K_n}$. If $d_n$ denotes the metric in $B_{K_n}$, then
$d_{n+1}|_{B_{K_n}}=d_n$, since $K_n\subset K_{n+1}$.
Then $\prod_{n=1}^{\infty }B_{K_n}=:Y$ in the Tychonoff product
topology is metrizable, that induces the metric in
$\Gamma _X$. Moreover, in the non-Archimedean case the metric $\rho $
in $Y$ can be chosen satisfying the ultrametric inequality:
$\rho (x,y):=d_n(x_n,y_n)p^{-n}$, where $n=n(x,y):=\min_{(x_j\ne y_j)}j$,
$x=(x_j: j\in {\bf N}, x_j\in B_{K_j} ) $.
\par As it was proved in the papers cited in \S 1, on $X$ from \S 2.1
there exist real measures $m $ quasi-invariant relative to
the left action of the corresponding group $G'$ such that $m(K_n)<\infty $. 
In the case $X=G$, then $G'$
is a dense subgroup in $G$. Quasi-invariance of $m $ implies, that
$m $ are non-atomic. Let $K\in \{ K_n: n\in {\bf N} \} $, then $m_K$
denotes the restriction $m|_K$. 
Then $m^n_K:=\bigotimes_{j=1}^n(m_K)_j$ is a measure
on $K^n$ and hence on ${\tilde K}^n$, since $m$ are non-atomic, where 
$(m_K)_j=m_K$ for each $j$. Therefore, $P_{K,m}:=exp(-m(K))\sum_{n=0}
^{\infty }m_{K,n}/n!$ is a probability measure on $Bf(B_K)$, where
$m_{K,0}$ is a probability measure on the singleton $B^0_K$, and
$m_{K,n}$ are images of $m^n_K$ under the following mappings:
$p^n_K: (x_1,...,x_n)\in {\tilde K}^n\to \{ x_1,...,x_n \} \in B^n_K$.
It was shown in \S 1.2 \cite{shim} that such system of measures
$P_{K,n}$ is consistent, that is, $\pi ^n_lP_{K_l,m}=P_{K_n,m}$ for each
$n\le l$. This defines the unique measure $P_m$ on $Bf(\Gamma _X)$, which is
called the Poisson measure.
For each $n_1,..,n_l\in \bf N_o$ and disjoint Borel subsets
$B_1,...,B_l$ in $X$ there is the following equality:
\par $(i)$ $P_m(\bigcap_{j=1}^l \{ \gamma : card(\gamma \cap B_i)=n_i \}
)=\prod_{i=1}^lm(B_i)^{n_i}exp(-m(B_i))/n_i!$.
\par The configuration space $\Gamma _X$ consists of
$\gamma \subset X$
such that $card(\gamma \cap K_n)<\aleph _0$ for each $n\in \bf N$. 
In the case of $Diff^t(M)$ this means that we need to consider
such elements $g$ of this group for which $supp(g)\subset K_n$ for some
$n\in \bf N$, for example, a subgroup with supports of its elements contained
in the corresponding finite unions of charts, where $supp(g):=cl \{ x\in M:
g(x)\ne x \} $. Such subgroups are not Banach
manifolds and they are denoted by $Diff_l(M)$. 
In the case of $X=G$ the initial configuration space
$\Gamma _X$ is not preserved by $G'$, since there are $g\in G'$ such that
$gK_n$ is not contained in any $K_m$, because $supp(L_h)=G$ for each $e\ne h
\in G'$, where $L_hg:=hg$ denotes the left shift in $G$ for $g, h \in G$.
\par Actually it is necessary to use more general construction
in the case of $X=G$. Let ${\tilde \Gamma }_X:=[\bigcup_{g\in G'}
g\Gamma _X]/R$,
where $R$ is an equivalence relation: $(g\gamma )R(g'\gamma ')$ 
if and only if $(g\gamma )=(g'\gamma ')$,
where $[\bigcup_{g\in G'} g\Gamma _X]$ is considered as the subset of
$X^{\bf N}$. The group $G'$ is separable,
hence there exists a countable dense subset $\{ g_j: j\in {\bf N} \} $.
To each element $g\in G'$ there corresponds a subsequence $\{ g_{j_n}:
n\in {\bf N} \} $ converging to $g$ in $G'$. Hence each $g\gamma $
is completely characterised by the corresponding subsequence 
$\{ g_{j_n}\gamma : n\in {\bf N} \} $. Therefore, ${\tilde \Gamma }_X$
has the embedding into $X^{\bf N}$ as the closed subset, 
since the family of mappings $\{ L_{g_j}: j\in {\bf N} \} $ separates
points of $X^{\bf N}$ \cite{eng}. Hence ${\tilde \Gamma }_X$
is also metrizable and complete. The manifold ${\tilde \Gamma }_X$
for each its point has a neighbourhood diffeomorphic with the corresponding
open subset of $\Gamma _X$, since for each $K_n$ there exist
a neighbourhood ${U'}_n$ of $e$ in $G'$ and $m>n$ 
such that ${U'}_nK_n\subset K_m$. A choice of such sequence $K_n\subset
Int(K_{n+1})$ with canonical closed subsets $K_n$ is given 
independently in \S 2.9.
The manifold ${\tilde \Gamma }_X$
is paracompact, consequently, it has a locally finite covering
$\{ S_l: l\in {\bf N} \} $, where $S_l$ are open in ${\tilde \Gamma }_X$
and diffeomorphic with the corresponding open subsets $Q_l$ of $\Gamma _X$
for which $P_m(Q_l)<\infty $, $\zeta _l: S_l\to Q_l$ denote such 
diffeomorphisms.
This means that the Poisson measure 
$P_m$ on $\Gamma _X$ induces the corresponding $\sigma $-additive
$\sigma $-finite quasi-invariant relative to $G'$ measure $\mu $ on ${\tilde
\Gamma }_X$ such that $\mu (E):=C\sum_lP_m(\zeta _l(E\cap S_l))2^{-l}$ 
which is also denoted by $P_m$, where $E\in Bf({\tilde \Gamma }_X)$,
a constant $C>0$ is chosen such that $\mu ({\tilde \Gamma }_X)=1$.
Therefore, $P_m$ on ${\tilde \Gamma }_X$ is the probability measure
as also for the case $\Gamma _M$ for $dim_{\bf L}M<\infty $.
This gives possibility to consider the case $X=G$ as well as the case
$X=M$ for groups of diffeomorphisms $Diff^t(M)$ of class $C^t$ with
$1\le t\le \infty $, 
which have structure of Banach manifolds from the
papers cited above.
\par If the manifold $M$ is locally compact and each $K_n$
is chosen to be canonical closed compact subset, then for $Diff^t(M)$
the configuration spaces $\Gamma _M$ and ${\tilde \Gamma }_M$
coincide. Indeed, if $\gamma \in \Gamma _M$, then $card(\gamma \cap K_n)
<\aleph _0$ for each $n\in \bf N$. Each subset $K_n$ is compact and canonically
closed, hence is sequentially compact \cite{eng}. This means that if
$card ((g\gamma )\cap K_l)=\aleph _0$ for some $l\in \bf N$ and 
$g\in Diff^t(M)$,
then $\{ g\gamma _j: j\in {\bf N} \} $ contains a convergent subsequence in
$K_l$. But $\{ \gamma _j: j\in {\bf N} \}=\gamma $ 
is the disrete subset of $M$,
hence $g^{-1}$ is not continuous, since $\{ g\gamma _j: j\in {\bf N} \}$
is not closed in $M$. This contradicts supposition $g\in Diff^t(M)$, consequently,
$g\gamma \in \Gamma _M$ for each $g\in Diff^t(M)$ for locally compact $M$ and 
canonical closed compact subsets $K_n$ in $M$. Therefore, $\bigcup_{g\in
Diff^t(M)}g \Gamma _M=\Gamma _M$, since $g\Gamma _M\subset \Gamma _M$ for each
$g\in Diff^t(M)$ and $e\Gamma _M=\Gamma _M$, consequently, 
${\tilde \Gamma }_M =\Gamma _M$.
\par If $M$ is not locally compact, for example, $M\setminus
M_R=\bigcup_{j=1}^{\infty }\Omega _j$, where $\Omega _j$ are disjoint open
subsets of $M$, $M_R:=\{ x\in M: d_M(x,x_0)\le R \} $, $0<R<\infty $,
$x_0$ is a fixed point in $M$ and $d_M$ is a metric in $M$,
then there exists $g\in Diff^{\infty }(M)$ with $supp(g)$ bounded in $M$
and a bounded infinite sequence of $\gamma _j\in M\setminus M_R$
which is discrete in $M$, that is, $cl \{ \gamma _j: j\in {\bf N} \}=
\{ \gamma _j: j\in {\bf N} \} $, such that $card ((g\gamma )\cap K_n)=
\aleph _0$ for some canonically closed $K_n$ in $M$, since each $K_n$ is not
locally compact, when $dim_{\bf L}M=\infty $, where $\bf L$ is the corresponding
field either $\bf R$ or the local field.
Hence in this case ${\tilde \Gamma }_M\ne \Gamma _M$.
\par If $X=G$, then in view of the choice of $K_n$ in \S 2.9 that to fulfil
demands on the measure $m$, there exists $g\in G'$ and $n\in \bf N$ such that
$gK_n$ is not contained in each $K_l$, where $l\in \bf N$. This $g$ can be
chosen by induction, since $K_l$ are not locally compact for each $l$ and $G$
is not locally compact. Therefore, there exists a discrete
infinite sequence $\gamma $ in $gK_n$ such that $card(\gamma \cap K_l)
<\aleph _0$ for each $l\in \bf N$. But $\gamma \in \Gamma _G$ and
$g^{-1}\in G'$ and $g^{-1}\gamma \in {\tilde \Gamma }_G\setminus \Gamma _G$,
since $card((g^{-1}\gamma )\cap K_n)=\aleph _0$. Hence ${\tilde \Gamma }_G
\ne \Gamma _G$ in this case also.
\par The group $G'$ and $X$ and $\Gamma _X$ have structures 
of the $C^{\infty }$-manifolds, since $X$ is the $C^{\infty }$-manifold.
Therefore, ${\tilde \Gamma }_X$ is the $C^{\infty }$-manifold also.
In the non-Archimedean case $M$, $G'$, $G$ and hence $\Gamma _X$ and
${\tilde \Gamma }_X$ are analytic manifolds with disjoint clopen charts,
since ${\Gamma '}_{\bf L}$ is discrete in $(0,\infty )$ and $\Gamma _X$ 
and ${\tilde \Gamma }_X$ are infinite-dimensional over $\bf L$
\cite{luum985}.
\par It is necessary to note, that for $X=G$ the dense subgroup $G'$
acts by the left shifts $L_h: G\to G$ as the diffeomorphism for each
$h \in G'$, where $G$ is either the loop group or the diffeomorphisms group. 
Therefore, lemmas 2.1, 2.2 and theorems 2.1, 2.2 and 2.3 from \cite{shim} are 
applicable to the cases considered here, since $\Gamma _X$ produces
charts for ${\tilde \Gamma }_X$ and $P_m$ on $\Gamma _X$ induces
$P_m$ on ${\tilde \Gamma }_X$. 
Theorem 2.3 from \cite{shim} can be applied to the real and non-Archimedean
cases of $X=M$.
\par {\bf 2.5. Definition. (see \cite{fell} \S 19.3.)} Let $G'$ 
be a group acting from the left
on the measure space $<X,{\sf L},m>$. Then $<X,{\sf L},m>$ is called
a measure $G'$-transformation space if $(i)$ $xW\in \sf L$ whenever
$x\in G'$ and $W\in \sf L$, and $(ii)$ $m(xW)=0$ whenever $x\in G'$,
$W\in \sf L$ and $m(W)=0$.
\par {\bf 2.6. Note.} For the considered here cases and $Bf(X)\subset \sf L$
conditions of definition 2.5 are fulfilled for the quasi-invariant measure
$m$ on $X$ relative to $G'$.
\par {\bf 2.7. Definition.} The measure $G'$-transformation space
$<S,{\sf L},m>$ is ergodic under $G'$ if, whenever, $V, W\in \sf L$
with $m(V)m(W)\ne 0$, there exists $x\in G'$ such that
$m(xV\cap W)\ne 0$.
\par {\bf 2.8. Note.} It was proved in
\cite{lutmf,luseamb,lurimut,luihlgr,luihlgna,luihgdsw}
that $m$ on $X$ is ergodic under $G'$ for the considered here cases
$(1-3)$, since $m$ is quasi-invariant relative to $G'$.
In cases $(1,2)$ at first $m$ 
was constructed on a neighbourhood $W$ of $e$ in $G$.
But theorem 2.3 from \cite{shim}
can not be applied to the cases $X=G$ for the probability measure $m$
on $X$, since in view of the construction
of the Gaussian measure $m$ on $G$ there are $\epsilon >0$ 
and $n\in \bf N$ such that for each $\psi \in G'$
with $\psi (K_n)\cap K_n=\emptyset $ the following integral is 
rather large: $\int_G
|\rho ^{1/2}_m(\psi ,x)-1|m(dx)>\epsilon ,$
where $m^{\psi }(E):=m(\psi ^{-1}E)$ for each $E\in Af(X,m)$,
$\rho _m(\psi ,x):=m^{\psi }(dx)/m(dx)$.
\par There are locally finite coverings $\{ g_jW_j: j\in {\bf N_o} \} $
of $G$ and $\{ g_jW'_j: j\in {\bf N_o} \} $ of $G'$, since $G$ and $G'$
are paracompact spaces relative to their own topologies
$\tau $ and $\tau '$ respectively and $G'$ is dense in $G$, where
$W_0=W$, $W_j\subset W$ for each $j$, $W'\subset W\cap G'$,
$W'_0=W'$, $W'_j\subset W'$ for each $j$, $g_j\in G'$ for each $j$,
$g_0=e$, $W_j$ are open in $G$ and $W'_j$ are open in $G'$.
Analogously for the pair $G'$ and $X=M$ in cases $(3)(a,b)$.
Then $m$ on $W$ can be extended as a $\sigma $-finite measure on $Bf(G)$
by the formula: 
\par $(i)$ $m(V):=\sum_{j=0}^{\infty }m(g_j^{-1}(V\cap g_jW_j))$,
since $0<m(W)<\infty $. The group $G$ is not locally compact,
hence $m(G)=\infty $. Using analogous procedure 
with a locally finite covering $\{ g_jW_j: j\in {\bf N_o} \} $
with $W_j$ open in $M$ and a neighbourhood $W$ of a marked point
$x_0\in M$ without relation between $W'$ and $W$
we get a $\sigma $-finite 
measure $m$ on $M$ for non-locally compact manifold $M$ with $m(M)=\infty $.
We choose in these cases $m(K_n)<\infty $ for each $n\in \bf N$.
As follows from the cited papers it is possible to choose $K_n
\subset \bigcup_{j=0}^ng_jW_j$ and $m$
such that 
\par $(ii)$ for each $\epsilon >0$ and each $n\in \bf N$ there exists
$\psi \in G'$ such that $\psi (K_n)\cap K_n=\emptyset $
and $\bigcup_nK_n=X$ and $\int_X|\rho ^{1/2}_m(\psi ,x)-1|^2m(dx)<\epsilon $.
Then it is proved below in theorem 2.9 that such $m$ exists and
$P_m$ on ${\tilde \Gamma }_X$ is ergodic.
Henceforth, such $\sigma $-finite measures $m$ on $X$ are used
with $m(X)=\infty $, since for $m(X)=1$ the corresponding measures $P_m$
are not ergodic (see note after definiton 1 in \S 2 \cite{shim}).
\par {\bf 2.9. Theorem.} {\it There exist 
quasi-invariant $\sigma $-finite measures $m$ on $X$ relative to
the groups $G'$ with $m(X)=\infty $ satisfying condition $(ii)$ 
from \S 2.8. For such $m$
the Poisson measure $P_m$ on ${\tilde \Gamma }_X$ is ergodic .}
\par {\bf Proof.} To prove $P_m$ is ergodic on $\Gamma _X$
we use the fact, that $m$ is ergodic on $X$. 
The measure space $<S,{\sf L},m>$ is said to have property $(P)$
if, for any locally $m$-measurable subset $W$ of $S$ such that
$xW\ominus W$ is locally $m$-null for each $x\in G'$, either
$W$ is locally $m$-null or $S\setminus W$ is locally $m$-null.
The measure space $<S,{\sf L},m>$ is called parabounded if there exists
a pairwise disjoint subfamily $\sf W$ of $\sf L$ such that
$(i)$ for each $A\in \sf L$, $\{ B\in {\sf W}: A\cap B\ne \emptyset \} $
is countable, and $(ii)$ $X\setminus \bigcup_{W\in \sf W}W$
is locally $m$-null.
It was proved in
proposition 19.5 \cite{fell} that if $<S,{\sf L},m>$
is ergodic it has property $(P)$. Conversely, if $<S,{\sf L},m>$
has property $(P)$ and is parabounded, it is ergodic.
The space $\Gamma _X$ is isomorphic with the projective limit
$pr-\lim \{ B_{K_n}, \pi ^n_m, {\bf N} \} $, which is the closed
subset in $\prod_nB_{K_n}$. The latter is the Polish space,
hence ${\tilde \Gamma }_X$ is the Polish space \cite{eng}.
The measure spaces $<X,Af(X,m),m>$ and $<{\tilde \Gamma }_X,
Af({\tilde \Gamma }_X,P_m),P_m>$  
are parabounded, since $X$ and ${\tilde \Gamma }_X$ are the Polish spaces and 
hence are the Radon 
spaces (see chapter 1 in \cite{dal}), that is, the class of compact subsets
approximates from below the corresponding measures $m|_{K_n}$ and $P_m$.
Therefore, it remains to show, that $<{\tilde \Gamma }_X,Af(
{\tilde \Gamma }_X,P_m),P_m>$
has property $(P)$.
But this follows from theorem 2.3 \cite{shim} and \S 2.8,
if to show that condition 2.8 (ii) is fulfilled for $m$
and ${U'}_nK_n\subset Int(K_{n+1})$ for the corresponding $K_n$ in $X$
and neighbourhoods ${U'}_n$ of $e$ in $G'$,
since there are the local diffeomorphisms $\zeta _l: S_l\to Q_l$
from \S 2.4 and $P_m$ and $m$ are $\sigma $-finite measures.
In this situation integral equalities and inequalities from the proof
of theorem 2.3 \cite{shim} are transferable onto the case of 
${\tilde \Gamma }_X$ considered here.
\par For the construction of such $m$ take it at first on an open subset
$U\subset X$ such that $W$ is sufficiently small: $W'W
\subset U$. In the case of $G=X$ in addition let $e\in U$
and $U^{-1}=U$, $W^{-1}=W$, $W'^{-1}=W'$ (see references in \S 2.8). 
The quasi-invariance 
factor $\rho _m(x,y)$ is continuos on $W'\times W$ and $\rho _m(e,y)=1$,
where $\rho _m(x,y):=m^x(dy)/m(dy)$, $m^x(A):=m(x^{-1}A)$
for each $x\in G'$ and $A\in Af(W,m)$.
Take open subsets $W_0\subset W$ and $W'_0\subset W'$ for which
$|\rho _m(x,y)-1|<1$ for each $(x,y)\in W'_0\times W_0$.
\par The measure $m$ is regular and approximated from above by the class
of open subsets \cite{dal,fell}. 
Therefore, it is possible to choose by induction
open subsets $W_j\subset W_0$ and $e\in W'_j\subset W'_0$ and
a sequence of elements $g_j\in G'$
such that $m(g_j^{-1}(g_jW_j\cap [\bigcup_{i=1}^{j-1}g_iW_i]))
<2^{-j}$ and $|\rho _w(x,y)-1|<2^{-j}$ on $g_jW'_j\times g_jW_j$,
where $w$ is a measure on $Bf(g_jW_j)$
defined by the following formula $w(g_jA):=m(A)$ 
for each $A\in Bf(W_j)$, $g_0=e$. Then $m$ on $Bf(X)$
is defined by formula 2.8(i) and certainly has the extension $m$ onto
$Af(X,m)$. 
\par The measure $m$ is induced from the corresponding measure
$\lambda $ on the Banach space $Y$
due to the local diffeomorphism $A: U\to V$, where $V$ is an open neighbourhood
of $0$ in $Y$ and $U$ is open in $X$. From the quasi-invariance of $\lambda $ 
relative to shifts from a dense subspace $Y'$ it follows a property:
\par $(\alpha )$ for each Borel subset $E\subset Y$ which is 
a $C^1$-submanifold in $Y$
of codimension $1$ in $Y$ (over the field
$\bf R$ or the non-Archimedean local field) such that
$T_yE$ is not subset of $Y'$ for each $y\in E$
it follows that $\lambda (E)=0$, since $\lambda $ is the quasi-invariant 
non-negative $\sigma $-additive and $\sigma $-finite measure.
In particular, for finite-dimensional $X=M$ over the corresponding field
the space $Y$ is finite-dimensional and $\lambda $ can be taken as
the Haar measure on $Y$ (in the real case it concides with the
Lebesgue measure). For infinite-dimensional real $X$, particularly for
$X=G$, the measure $\lambda $ can be taken Gaussian. For infinite-dimensional
$X$ over the local field the wider class of measures $\lambda $
was constructed in the papers cited in \S 1.
Then we choose (take) by induction a sequence $K_n\subset
\bigcup_{i=0}^ng_iW_i$ satisfying the following conditions
${U'}_nK_n\subset Int(K_{n+1})$ for each $n$
with $\bigcup_nK_n=X$ and $m(K_n\setminus Int(K_n))=0$
and $K_n$ are canonical closed subsets, that is, $cl(Int(K_n))=K_n$,
since $m$ is quasi-invariant and has not any atoms
and due to property $(\alpha )$ of $\lambda $,
where $cl(A)$ denotes the closure of a subset $A\subset X$ in $X$,
$Int(A)$ denotes the interior of $A$ in $X$, ${U'}_n$ are the corresponding
(open) neighbourhoods of $e$ in $G'$ such that ${U'}_n\subset W'$. 
The space $X$ is Polish,
hence each $K_n$ is the Polish topological subspace \cite{eng}.
Certainly, in the non-Archimedean cases each $K_n$ can be chosen clopen
(closed and open) in
$X$, that is, $Int(K_n)=K_n=cl(K_n)$, since the base of the topology of $X$
consists of clopen subsets.
Since $X$ is not locally compact, then there exists 
the sequence $\{ K_n: n\in {\bf N} \}$ fulfilling condition 2.8(ii).
\par {\bf 2.10. Note.} In cases $(3)(a,b)$ for $X=M$ and $G'=Diff^t(M)$ 
in addition we have the following.
\par {\bf 2.11. Definition.} Let ${G'}_{K_n}:=\{ \psi \in G':
\psi |_{K^c_n}=id \} $ and let $f$ be a symmetric measurable function 
defined on ${\tilde K}^l_n$, where $l\in \bf N$, $A^c:=X\setminus A$
for a subset $A$ in $X$, $K_n$ are canonical closed subsets
with $\bigcup_nK_n=X$ and $K_n\subset K_{n+1}$ and $m(K_n\setminus Int(K_n))=0$
for each $n\in \bf N$. In the non-Archimedean case let also $K_n$ be clopen in
$M$, which automatically implies $K_n\setminus Int(K_n)=\emptyset $.
The measure $m$ is called ${G'}^l_{K_n}$-ergodic,
if $f$ is constant modulo null sets, then $f(x_1,...,x_l)=f(\psi (x_1),...,\psi
(x_l))$ for $m^l_{K_n}$-a.e $x=(x_1,...,x_l)$ for all $\psi \in G'$.
\par {\bf 2.12. Theorem.} {\it If for each $n$ the measure $m$ 
is ${G'}^l_{K_n}$-ergodic for some $N\ge n$ and all $l$, then $P_m$
is $G'$-ergodic.}
\par {\bf Proof.} As it was shown in papers
\cite{lutmf,lurimut,lupr202,lupr215,lupr218} the subgroups ${G'}_{K_n}$
are correctly defined for canonical closed subsets $K_n$ in $M=X$,
${G'}_{K_n}\subset G'$ for each $n$, since from $\psi |_{K^c_n}=id$
it follows, that $\psi |_{cl(K^c_n)}=id$.
The rest of the proof is as in the proof of theorem 2.4 \cite{shim},
which can be applied locally and then with the help of the local
diffeomorphisms $\zeta _l: S_l\to Q_l$ is extendable onto the case of
${\tilde \Gamma }_X$ considered here, since ${G'}^l_{K_N}{\tilde K}^l_N=
{\tilde K}^l_N$ and for the measure $\nu (A):=P_m(E\cap A)$ for each
$A\in Bf({\tilde \Gamma }_X)$ we have $\nu (B)=\int_0^{\infty }P_{cm}(B)
\lambda (dc)$ for each $B\in Bf(\Gamma _X)$,
where $c\ge 0$ and $\lambda $ is a suitable Borel measure on $[0,\infty )$.
From $P_m(\zeta _l(A\cap S_l))=0$ for each $l$ it follows, that
$P_m(A)=0$. Thus if $\lambda (\{ 1 \} )>0$, then $P_m(A)=0$;
if $\lambda (\{ 1 \} )=0$, then $P_m(A^c)=0$, where $A$ is a $P_m$-measurable 
subset of ${\tilde \Gamma }_X$ for which $P_m(A\bigtriangleup \psi ^{-1}A)=0$ 
for all $\psi \in G'$, where $A\bigtriangleup B:=(A\setminus B)\cup (B\setminus
A))$.
\par {\bf 2.13. Note.} From theorem 2.12 it can be deduced in another way, 
than it was done in theorem 2.9, that $P_m$ on ${\tilde \Gamma }_X$
is $G'$-ergodic in cases $(3)(a,b)$
for $X=M$, when $m$ and $K_n$ are chosen in accordance with \S 2.8 and \S 2.11.
The proof of this is analogous to that of theorem 2.5 \cite{shim},
since $m$ is ergodic and quasi-invariant with the continuous quasi-invariance
factor $\rho _m(\psi ,x)$ on $G'\times X$, 
$m(X)=\infty $ and $m(K_n\setminus Int(K_n))=0$, since due to \S 2.4
there are the local diffeomorphisms $\zeta _l: S_l\to Q_l$
and $Diff^t(X;K){\tilde K}^l={\tilde K}^l$ for each canonical closed subset $K$
in $X$, where $Diff^t(X;K):=\{ f\in Diff^t(X): f|_{K^c}=id \} $.
\par {\bf 2.14. Lemma.} {\it Let $Y$ be a canonically closed subset in $X$,
$Y\subset K_n$ for some $n\in \bf N$. Suppose that $\mu $ is a quasi-invariant
measure on ${\tilde \Gamma }_X$ relative to $G'=Diff^t(X)$ for a $C^{\infty
}$-manifold $X=M$ (in the non-Archimedean case an analytic manifold $M$)
and $\mu _n$ be a restriction of $\mu $
on $B^n_Y\times {\tilde \Gamma }_{X\setminus Y}$ 
and $\mu '_n$ and $\mu "_n$ be
projections of $\mu _n$ on $B^n_Y$ and ${\tilde \Gamma }_{X\setminus Y}$ 
respectively.
Then $\mu _n$ is equivalent with $\mu '_n\times \mu "_n$.
In the non-Archimedean case this is also true for $Y$ clopen in $X$.}
\par {\bf Proof.} In view of \S 2.9 $m(Y\setminus Int(Y))=0$. The group
$Diff^t(X;X\setminus Y)$ is a closed subgroup of $Diff^t(X)$, hence
${\tilde \Gamma }_{X\setminus Y}=(Diff^t(X;X\setminus Y)\Gamma _X)/R$
is a $C^{\infty }$-submanifold of ${\tilde \Gamma }_X$ (see also \S 2.4).
The measures $\mu '_n\times \mu "_n$ with $\mu _n$
are equivalent if and only if $\mu $ and $\mu '\times \mu "$ are equivalent,
since $\mu $ is quasi-invariant relative to $G'$ and non-atomic, where $\mu '$
is a projection of $\mu |_{B_Y\times {\tilde \Gamma }_{X\setminus Y}}$ on 
$B_Y$ and $\mu "$ is a projection of $\mu $ on 
${\tilde \Gamma }_{X\setminus Y}$. 
On the other hand, $G'B^n_X=B^n_X$ for each $n\in \bf N$
and $G'{\tilde \Gamma }_X={\tilde \Gamma }_X$, also $Diff^t(X;X\setminus Y)
{\tilde \Gamma }_{X\setminus Y}={\tilde \Gamma }_{X\setminus Y}$.
On the other hand, $B^n_Y\times {\tilde \Gamma }_{X\setminus Y}$ is the
Borel subset of ${\tilde \Gamma }_X$, since $B^n_Y$ is the Borel subset
of $\Gamma _X$. For the rest of the proof are necessary two propositions.
\par {\bf 2.15. Proposition.} {\it In the group 
$Diff^t(Int(Y))$ there exists a countable
family of one-parameter subgroups $G_l$ such that generated by them
group $J\subset Diff^t(Y)$ acts transitively on $B^n_Y$.}
\par {\bf Proof.} For $Diff^t(Int(Y))$ one-parameter 
subgroups can be chosen as in proposition 2.1 \cite{vgg} with the help
of \cite{ebi} and theorems about existence of one-parameter subgroups
of $Diff^t(Y)$ for infinite-dimensional Banach manifolds $M$ 
from \cite{luumn96,lurimut,lupr180,lupr202}, where one-parameter subgroups are
real for real $M$ and $g^b$ with $b\in \bf L$ for $M$ over the local field
$\bf L$ such that $g^ag^b=g^{a+b}$ for each $a, b\in \bf L$.
In the non-Archimedean case one-parameter subgroups can aslo be indexed
by $b\in B({\bf L},0,1)$, where $B(S,x,r):=\{ y\in S: d_S(x,y)\le r \} $
denotes a ball in 
a metric space $S$ with a metric $d_S$ and a point $x\in S$. 
This is possible, since $M$ and $T_xM$ are separable spaces for each $x\in M$
and using countable atlas $At(M)=\{(U_j,\phi _j): j \} $ of $M$
and considering one-parameter subgroups with $supp(g^b)\subset U_j$
for each $b\in \bf L$ for the corresponding chart $U_j$, 
where either ${\bf L}=\bf R$ or $\bf L$ is the local
field, $U_j$ are open in $M$ and $\phi _j: U_j\to V_j$ are diffeomorphisms,
$V_j$ are open in the corresponding Banach space.
\par {\bf 2.16. Proposition.} {\it Let $\bf L$ or may be $B({\bf L},0,1)$
also in the non-Archimedean case 
acts measurably in a measure 
space $(M,Bf(M),\mu )$ such that $\mu $ is quasi-invarint relative to the 
action of $\bf L$ or $B({\bf L},0,1)$ on $M$, 
where $M$ is a $C^{\infty }$-manifold
(analytic in the non-Archimedean case) and $\mu $ is 
induced by a quasi-invariant non-negative
$\sigma $-additive and $\sigma $-finite measure $\eta $ relative to
shifts from a dense subspace $Z'$ and $\eta $ is on the Borel field $Bf(Z)$ 
of the separable Banach space $Z$ 
over a field $\bf L$ which is either ${\bf L}=\bf R$
or a local field such that $Z=T_xM$ for each $x\in M$. 
Suppose that a partition $\zeta $ of $M$ is invariant
by $mod (\mu )$ relative to the action of $\bf L$ or $B({\bf L},0,1)$
on $M$ and projections of
$\eta $ onto one-dimensional over $\bf L$ subspaces are equivalent with the 
non-negative Haar
measure $\lambda $ on $\bf L$. Then for 
$\mu $-almost each $C\in \zeta $ the conditional measures $\mu ^C$ are
quasi-invariant relative to the action of $\bf L$ or $B({\bf L},0,1)$
respectively.}
\par {\bf Proof.} The proof is almost the same as in proposition 2.2
\cite{vgg} with the substitution of $\bf R$ onto $\bf L$ 
or may be $B({\bf L},0,1)$ in the non-Archimedean case and using the Haar
measure $\lambda $ on a locally compact subgroup $S$ of $\bf L$ with
$\lambda ({\bf L}\setminus S)=0$ or $\lambda (B({\bf L},0,1)\setminus S)=0$,
which implies $S=\bf L$ or $S=B({\bf L},0,1)$ respectively by the A. Weil
theorem, since each quasi-invariant measure on a locally compact group
(relative to its action on itself) is equivalent with the Haar measure
\cite{boui}.
\par {\bf Continuation of the proof of lemma 2.14.} In view of proposition
2.15 there exists a subgroup $J$ which acts transitively on $B^n_{Int(Y)}$.
In view of proposition 2.16 from an isomorphism of one-parameter subgroup
$G_l$ with $\bf L$ or $B({\bf L},0,1)$ for $\mu "_n$-a.e. configurations
$\gamma \in {\tilde \Gamma }_{X\setminus Y}$ the conditional measure
$\mu ^{\gamma }_n$ on $B^n_Y$ is quasi-invariant relative to each
one-parameter subgroup $G_l$, hence realtive to the minimal subgroup
$J$ of $Diff^t(X)$ generated by $\bigcup_{l=1}^{\infty }G_l$.
The measure $\mu ^{\gamma }_n$ on $B^n_Y$ induces the measure $\eta $
on $T_{\gamma ^n}
B^n_Y$ for each $\gamma ^n\in B^n_Y$. This measure $\eta $ is completely
characterised by its finite-dimensional projections $\eta _n$ 
onto subspaces $F_n$
such that $F_n\subset F_{n+1}$ for each $n\in \bf N$ and $\bigcup_nF_n$
is dense in $T_{\gamma ^n}B^n_Y$ (see about weak distributions
\cite{lupr210,sko}). It is supposed that the manifold $M$ has the foliated
structure such that $M_n\subset M_{n+1}$ and
$dim_{\bf L}M_n=k_n<\infty $ for each $n\in \bf N$ and
$\bigcup_{n\in \bf N}M_n$ is dense in $M$. Theorefore, to $\mu ^{\gamma }_n$
there corresponds a family of measures $\tilde \eta _n$ on $M_n$
with the help of a locally finite coverings and
the exponential mapping $exp: \tilde M\to M$
from the neighbourhood $\tilde TM$ of $M$ in $TM$ onto $M$
such that $exp_x: V_x\to W_x$ are local diffeomorphisms of open subsets
$V_x$ in $T_xM$ and $W_x$ in $M$ with $x\in M$. A measure
$\tilde \eta _n$ is quasi-invariant relative to $J_n:=\{ g\in J: g_{
M\setminus M_n}=id \} $. The manifold $M_n$ is locally compact,
hence $\tilde \eta _n$ is equivalent with the Riemann volume element
on $M_n$ in the real case and with the restriction of the Haar
measure from ${\bf L}^{k_n}$ onto $M_n$ in the non-Archimedean case, since
in the latter case $M_n$ is embeddable into ${\bf L}^{k_n}$ due to a
partition of $M_n$ into a disjoint union of balls. In view of the Kakutani
theorem II.4.1 \cite{dal} $\mu _n$ is equivalent with $\mu '_n\times \mu "_n$,
since a finite measure $\zeta $ on $Bf(A\times B)$ 
is equivalent with the direct product
$\zeta _A\times \zeta _B$, where $\zeta _A$ and $\zeta _B$ are projections of 
$\zeta $ on Hausdorff topological spaces $A$ and $B$ respectively.
\section{Unitary representations associated with the Poisson measures.}
\par {\bf 3.1. Definitions and notes.} Let $H:=L^2({\tilde \Gamma }_X,
P_m,{\bf C})$ 
be the standard Hilbert space of equivalence classes of measurable
functions $f: {\tilde \Gamma }_X\to \bf C$ 
for which $\| f \|_H^2:=\int_{{\tilde \Gamma }_X}
|f(x)|^2P_m(dx)<\infty $, where $P_m$ is the Poisson measure given in \S 2.4.
Then consider the following representation:
\par $(i)$ $U_m(\psi )f(\gamma ):=\rho _{P_m}^{1/2}
(\psi ,\gamma )f(\psi ^{-1}(\gamma ))$,\\ 
where $\rho _{P_m}(\psi ,\gamma ):=P_{m^{\psi }}(d\gamma )/P_m(d\gamma )$,
$\gamma \in {\tilde \Gamma }_X$, $f\in H$, $\psi \in G'$, 
$m^{\psi }(E):=m(\psi ^{-1}E)$ for each $E\in Af(X,m)$.
That is,
$U_m: G'\to U(H)$, where $U(H)$ is the unitary group of the Hilbert space
$H$. The topology of $U(H)$ is induced by the operator norm in the space
$L(H)$ of bounded linear operators $S: H\to H$, $d(A,B):=d(B^{-1}A,I):=
\| B^{-1}A-I \|_{L(H)}$ is the metric in $U(H)$, where $A, B\in U(H)$,
$I$ denotes the unit operator on $H$.
\par In cases $(3)(a,b)$ of $X=M$ and $G'=Diff^t(M)$
these representations can be generalised
with the help of the symmetric group $\Sigma _n$ 
representations in the following 
manner, where $\Sigma _n$ 
is the group of all (bijective) automorphisms $\sigma $
of the set $\{ 1,2,..,n \}$ with $n\in \bf N$ and $\Sigma ^{\infty }$ 
is the symmetric
group of $\bf N$ (that is, of all bijective mappings of $\bf N$).
Let $q: \Sigma _n\to U(W)$ be the unitary representation of $\Sigma _n$, where
$W$ is the Hilbert space, or $q: \Sigma ^{\infty }\to U(W)$. 
Then $s_n: B^n_X\to \tilde X^n$ or $s: {\tilde \Gamma }_X\to {\tilde X}^{\infty
}$ produces
a mapping $\sigma : G'\times B^n_X\to \Sigma _n$ 
or $s: G'\times {\tilde \Gamma }_X
\to \Sigma ^{\infty }$ by the formula 
$s_n(\psi ^{-1}(\gamma ))=\psi ^{-1}(s_n(\gamma ))\sigma (\psi ,\gamma )$
or $s(\psi ^{-1}(\gamma ))=\psi ^{-1}(s(\gamma ))\sigma (\psi ,\gamma )$,
where $s_n$ is a measurable cross-section of $p_n : \tilde X^n\to B^n_X$
and $s$ of $p: {\tilde X}^{\infty }\to {\tilde \Gamma }_X$
such that $p_n(x_1,...,x_n)=\{ x_1,...,x_n \} $ for $n\in \bf N$
or $p(x_1,x_2,...)=\{ x_1,x_2,... \} $, 
$(x_1,...,x_n)\sigma =(x_{\sigma (1)},...,x_{\sigma (n)})$
or $(x_1,x_2,...)\sigma =(x_{\sigma (1)},x_{\sigma (2)},... )$
respectively. Then 
with each pair $(n,q)$ or $(\infty ,q)$ 
is associated a unitary representation of $G'$
in $L^2(B^n_X,m^n,W)$ or in $L^2({\tilde \Gamma }_X,P_m,W)$ respectively 
such that 
\par $(ii)$ $V^q_m(\psi )f(\gamma ):=\rho _{m^n}^{1/2}
(\psi ,\gamma )q(\sigma (\psi
,\gamma ))f(\psi ^{-1}(\gamma ))$, or\\
\par $(iii)$ $V^q_m(\psi )f(\gamma ):=\rho _{P_m}^{1/2}(\psi ,\gamma )
q(\sigma (\psi ,\gamma ))f(\psi ^{-1}(\gamma ))$,\\
where $m^n$ is the image measure of the direct product of $n$ copies of
$m$ by the map $p_n$ and $\rho _{m^n}(\psi ,\gamma ):=(m^n)^{\psi }(d\gamma )/
m^n(d\gamma )$, $(m^n)^{\psi }(E):=m^n(\psi ^{-1}E)$ for each 
$E\in Af(B^n_X,m^n)$, $\psi \in G'$, $\gamma \in B^n_X$.
As usually the space $L^2(B^n_X,m^n,W)$ denotes the space of equivalence
classes of measurable
functions $f: B^n_X\to W$ for which $\| f\| ^2_{L^2(B^n_X,m^n,W)}:=
\int_{B^n_X}\| f(x)\| _W^2m^n(dx)<\infty $.
Then $U_m$ and $V^q_m$ define new representation $U^q_m:=U_m\otimes V^q_m$.
\par {\bf 3.2. Note.} For the group of diffeomorphisms of the real
finite-dimensional manifold $M$ such representations were defined in
\cite{vgg}, where it was mentioned that the representations $V^q_m$
are in a weak respect
analogous to the construction of H. Weyl for the classical Lie groups.
For $W=\{ 0 \}$ and $q=I$ the representation $V^q_m$ is evidently irreducible
\cite{hira,lurimut,luseamb},
hence as in theorem 1.1 \cite{vgg} for the considered here cases
we have, that $(i)$ if $q$ is the irreducible representation of $\Sigma _n$
with $n\in \bf N$, then
$V^q_m$ is the irreducible representation of the diffeomorphism group $G'$;
$(ii)$ $V^{q_1}_m$ and $V^{q_2}_m$ are equivalent if and only if $n_1=n_2$
and $q_1$ of $S_{n_1}$ is equivalent to $q_2$ of $S_{n_2}$.
\par {\bf 3.3. Note.} Let $X=M$ be a finite-dimensional over a local field
$\bf K$ non-compact manifold embedded as an open subset into $\bf K^n$.
Suppose that $m$ is the restriction $m=\lambda |_M$ 
of the Haar measure $\lambda $ on $\bf K^n$ normalised by $\lambda (B({\bf
K^n},0,1))=1$. Let $Diff^t(X,m)$ denotes the subgroup of $G'=Diff^t(X)$ 
of the non-Archimedean class of smoothness $C^t$ such that 
$\rho _m(\psi ,x)=1$ for each $\psi \in Diff^t(X,m)$ and $x\in X$,
where $1\le t\le \infty $.
\par {\bf 3.4. Theorem.} {\it Let $X$ and $Diff^t(X,m)$ be the same as in \S 3.3.
Then the restriction of the representation
$V^q_m$ from \S 3.1 on $Diff^t(X,m)$ is irreducible.}
\par {\bf Proof.} For finite-dimensional $M$ over the local field
$\bf L$ there is the equality ${\tilde \Gamma }_M=\Gamma _M$ due to \S 2.4.
Since diffeomorphisms $\psi $ with locally linear 
$(\psi -id)$ are contained in $Diff^t(X,m)$, for example, when
$\| \psi -id \|_{C^1(X\to {\bf K^n})}<1$. Then for each pairwise distinct
points $x_1,...,x_n\in X$ there are neighbourhoods $O_1,...,O_n$
such that their closures $\bar O_j$ are $C^1$-diffeomorphic with 
balls in $\bf K^n$ and $\bar O_i\cap \bar O_j=\emptyset $ for each
$i\ne j$ and $m(O_1)=...=m(O_n)$. Moreover, for each transposition
$(k_1,...,k_n)$ of $(1,...,n)$ there exists a diffeomorphism $\psi
\in Diff^t(X,m)$ with $\psi (\bar O_i)=\bar O_{k_i}$. Such
$\psi $ exists due to partition of $X$ into disjoint union of sufficiently 
small clopen balls $U_l$ such that $\bar O_j=O_j$ for each $j=1,...,n$
and for each $j$ there exists $l_j$ such that $O_j=U_{l_j}$
while $\psi (O_j)=O_{k_j}$ and diameters of all $O_j$ are equal to each 
other. Having $\psi $ on $\bigcup_{j=1}^nO_j=:E$, it is possible to extend
$\psi $ as $id$ on $X\setminus E$. 
\par Let now $Y$ be a clopen compact 
submanifold of $X$ and consider subspace ${\tilde L}^2(Y,m,{\bf C})$
consisting of $f\in L^2(Y,m,{\bf C})$ with $\int_Yf(y)m(dy)=0$.
Let $H_1\ne \{ 0 \}$ be an invariant subspace of ${\tilde L}^2(Y,m,{\bf C})$
relative to the regular representation
$(U_{\psi }f)(y):=f(\psi ^{-1}y)$ of $Diff^t(Y,m)$.
For each ball $O$ in $Y$ there exists $f\in H_1$ such that $f\ne 0$ and
$supp(f)\subset O$. Further analogously to the proof of lemma 2 from \S 1 
\cite{vgg} we get, that such representation of $Diff^t(Y,m)$ is irreducible.
If $O_j$ with $j=1,...,n$ are clopen balls in $X$, then a subgroup
$G^0(O_1,...,O_n;X)$ of elements $\psi \in Diff^t(X,m)$ with $\psi |_{O_j}=id$
for each $j=1,...,n$ acts trivially on $(\bigotimes_{j=1}^nL^2(O_j,m,{\bf C})
)\otimes W$. Then quite analogously to the proof of theorem 1.2 \cite{vgg}
we get the statement of this theorem.
\par {\bf 3.5. Note.} Let ${\tilde {\bf N}}^n:=\{ a=(i_1,...,i_n)|
i_j\ne i_s$ for each $i\ne s \} $, $l_2({\tilde {\bf N}}^n,W):=
\{ \phi | \phi : {\tilde {\bf N}}^n\to W,$ such that $\| \phi \|^2:=
\sum_{a\in {\tilde {\bf N}}^n}\| \phi (a)\| ^2_W<\infty \} $
and $H^q:=\{ \phi \in l_2({\tilde {\bf N}^n},W)| \phi (i_{\sigma (1)},...,
i_{\sigma (n)})=$ $q^{-1}(\sigma )\phi (i_1,...,i_n))$ for all $\sigma \in
S_n \} $, where $q$ is a unitary representation of $S_n$ in a Hilbert space
$W$. In the case $X=M$ the representation $q$ and $W$ may be non-trivial
with $n\in \bf N$,
for $X=G$ and $G'$ acting on $X$ we set $q=I$ and $W=\{ 0 \}$ and $n=0$.
We denote by $\Sigma ^{\infty }$ the set of all permutations (bijections)
of $\bf N$ and put $\sigma a=(\sigma (i_1),...,\sigma (i_n))$ for
$\sigma \in \Sigma ^{\infty }$ and $a\in {\tilde {\bf N}^n}$.
Then a function $\sigma : G'\times {\tilde \Gamma }_X
\to \Sigma ^{\infty }$ is defined 
by the formula $s(\psi ^{-1}(\gamma ))=\psi ^{-1}(s(\gamma ))\sigma (\psi,
\gamma )$, where $s$ is a measurable (admissible) cross section
of the map $p: {\tilde X}^{\infty }\ni (x_1,x_2,...)\mapsto \{ x_1,x_2,... \}
\in {\tilde \Gamma }_X$ 
posessing the following property: 
\par $(\alpha )$ if $card (\gamma \cap
X_1)=k_1$, $card (\gamma \cap (K_2\setminus K_1))=k_2$,..., $card (\gamma 
\cap (K_n\setminus K_{n-1}))=k_n$,..., then the first $k_1$ elements of
$s(\gamma )$ are in $\gamma \cap K_1$, the next $k_2$ of $s(\gamma )$
are in $\gamma \cap (K_2\setminus K_1)$ and so on.
For $X=M$ finite-dimensional over $\bf L$
and $\psi \in Diff^t(K_l):=\{ \psi \in Diff(X): \psi |_{K^c_l}=id \} $,
then $\sigma (\psi ,\gamma )\in \Sigma _r$, so $H^q$ is non-trivial in general
for this case, where $K^c_l:=X\setminus K_l$. The latter property in general
may be untrue for infinite-dimensional manifold $M$ or for
$X=G$ and $G'$ acting on $G$, therefore, we consider
$q=I$ and $n=0$ and $W=\{ 0\} $ for $X=G$.
For infinite-dimensional $X=M$ over $\bf L$ let us drop condition $(\alpha )$
and let $q$ be a representation of $\Sigma ^{\infty }$ in $U(H^q)$,
where $H^q$ is defined analogously with the $H^q$ for $n$ but with the
substitution of $n$ onto $\infty $ and $\Sigma _n$ onto $\Sigma ^{\infty }$.
\par Then there exists the following unitary representation
of $G'$ in the space $L^2({\tilde \Gamma }_X,P_m,{\bf C})\otimes H^q$
(which is isomorphic with $L^2({\tilde \Gamma }_X,P_m,{\bf C})$ for $X=G$):
\par $(i)$ $U^q_m(\psi )F(\gamma ,a):=\rho _{P_m}^{1/2}(\psi ,\gamma )
F(\psi ^{-1}(\gamma ),\sigma (\psi ,\gamma )^{-1}a)$
(see $\rho _{P_m}$ in \S 3.1).
\par {\bf 3.6. Theorem.} {\it The representations from \S 3.1 and \S 3.5
in the case $X=G$ are equivalent, in the case of finite-dimensional
$X=M$ over $\bf L$ for the group of
diffeomorphism acting on $M$ the representation $U_m\otimes V^q_m$
is equivalent with $U^q_{n\circ m}$, where $q$ is a representation
of the symmetric group $\Sigma _n$.}
\par {\bf Proof.} In the case $X=G$ this follows from their definitions,
that is, $U_m$ given by formula $3.1.(i)$ is equivalent with
$U^q_m$, since $q=I$ and $W$ is trivial. In the case of the $C^{\infty
}$-manifold $X=M$, which is finite-dimensional over $\bf R$, it was proved in
theorem 3.2 \cite{vgg}.
In the non-Archimedean case the proof is analogous, but instead of
differentiabilty of measures their pseudodifferentiability
should be considered as in \cite{lupr210}. In view of proposition
2.16 the quasi-invariant measure $m$ on $M$ relative to $Diff^t(M)$
is equivalent with the restriction of the Haar measure $\lambda $
for ${\bf L}^k$ on $M$, that is, $\lambda |_M \sim m$,
where $M$ is embedded into ${\bf L}^k$
for the corresponding local field $\bf L$.
In view of \S 2.4 $\Gamma _M={\tilde \Gamma }_M$ for locally compact $M$.
Therefore, convolutions of measures $\mu _1 *\mu _2$ are correctly defined
on $\Gamma _M$ as an image of the product measure $\mu _1\times \mu _2$
on $\Gamma _M\times \Gamma _M$ 
relative to a mapping $(\gamma _1,\gamma _2)\mapsto (\gamma _1\cup \gamma _2)$
from $\Gamma _M\times \Gamma _M$ to $\Gamma _M$.
Then $n\circ \mu $ denotes $\mu *m_n$, where $m_n$ is a measure on $B^n_X$
corresponding to the restriction $P_m|_{B^n_X}$, $n_1\circ (n_2\circ \mu )$
is equivalent with $(n_1+n_2)\circ \mu $ for each $n_1$ and $n_2
\in \bf N$, $0\circ \mu $ is equivalent with $\mu $. 
For each $\psi \in Diff^t(M)$ we have $\psi (\mu _1*\mu _2)=(\psi \mu _1)*(\psi
\mu _2)$, where $t\ge 1$ and $\psi \mu :=\mu ^{\psi }$,
$\mu ^{\psi }(E):=\mu (\psi ^{-1}E)$ for each $E\in Bf(M)$. Therefore,
for each pair $\mu _1$ and $\mu _2$ of quasi-invariant measures, their
convolution is also a quasi-invaraint measure. Then all necessary results from
\S 2.3-5 of \cite{vgg} can lightly be
transferred onto the non-Archimedean case.
\par {\bf 3.7. Note.} In the papers \cite{luummr,lupr202,lupr218,shav}
quasi-invariant measures on the diffeomorphisms groups of real Banach
manifolds were constructed. Purely Gaussian measures 
quasi-invariant relative to dense subgroups were constructed
in the cases of Euclidean and Hilbertian at infinity manifolds and also
for definite closed subgroups $Diff^t_k(M)$ $:=\{ f\in Diff^t(M):
(\Delta ^jf)|_{\partial M}=(\Delta ^j)id|_{\partial M}$ $\mbox{for each}$ 
$j=0,1,...,k \} $ and $Dif^t_l(M):=\{ f\in Diff^t(M): 
(\partial _{\nu }^jf)|_{\partial
M}=(\partial _{\nu }^jid)|_{\partial M}$ $\mbox{for each}$ $j=0,1,...,l \} $
of $Diff^t(M)$ for compact $C^s$-manifolds $M$
with a boundary $\partial M\ne \emptyset $, where $Diff^t(M)$ has
a class of smoothness
by H\"older $C^t$, also a class of smoothness $H^t$ by Sobolev or Besov
was considered for $t>dim_{\bf R}M+5$, $\Delta $ denotes the
Beltrami-Laplace operator on $M$, $\partial _{\nu }$ denotes the partial
differentiation along normal to the boundary in local coordinates,
$\Delta ^0=I$ and $\partial ^0=I$ are the unit operators. 
In particular for a compact manifold with
the boundary purely Gaussian measures $\mu $ on $Diff^t_k(M)$ and $Dif^t_l(M)$
quasi-invariant relative to dense subgroups $Diff^{t'}_k(M)$
and $Dif^{t"}_l(M)$ were constructed, where $k$ and $l$ and $t'-t>0$
and $t"-t>0$ are dependent on $dim_{\bf R}M$, $s>t'+2$ and $s>t"+2$
respectively. The cases of Schwarz
class of smoothness also were considered. The given below theorem
in the real case for finite measures
was proved shorlty earlier in \cite{luihgdsw},
in the non-Archimedean case it is contained in \cite{lutmf}.
In theorem 3.8 a quasi-invariant $\sigma $-finite $\sigma $-additive
measure is considered, which may be unbounded.
The cases of $\sigma $-finite non-negative 
measures and probability measures on
$G$ are considered quite analogously. Certainly this theorem is applicable not
only to Gaussian measures but also to measures which have definite properties of
the quasi-invariance factors $\rho _{\mu }$ such that 
a family of continuos 
functions $\{ \rho _{\mu }^{1/2}(z,g)=\phi (g): z\in G' \} $
parametrized with $z\in G'$ separates points of $G$ (see more precisely
the proof below). It is essential in the proof
that $G$ is the infinite-dimensional 
non-locally group and $G'$ is its dense subgroup such that the measure
$\mu $ is ergodic. Evidently, if $\mu '$ is a measure equivalent with
$\mu $, then the regular representations associated with them are equivalent
due to the isomorphism $f(g)\mapsto (\mu '(dg)/\mu (dg))^{1/2}f(g)$
of the Hilbert space $L^2(G,\mu ',{\bf C})$ with $L^2(G,\mu ,{\bf C})$,
where $f\in L^2(G,\mu ',{\bf C})$ and $g\in G$.
\par {\bf 3.8. Theorem.} {\it Let $G$ be a group 
of diffeomorphisms with a real probability quasi-invariant
measure $\mu $ relative to a dense subgroup $G'$ as in \S 3.7.
Then $\mu $
may be chosen such that the associated regular unitary representation 
of $G'$ is irreducible.}  
\par {\bf Proof.}  Let a measure $\nu $ on a Banach space
$H$ be of the same type as in the proofs of theorems in papers cited in
\S 3.7 such that a local diffeomorphism $A: W\to V_H$ induces a quasi-invaraint
measure on $W$ and then with the help of left shifts $g_j\in G'$ on the entire 
group $G$, where $W$ is an open neighbourhood of $e$ in $G$ and $V_H$ is an open
neighbourhood of $0$ in $H$. We choose a constant multiplier $c>0$ for $\mu $
such that $c\mu (W)=1$ and then denote such normalized measure by $\mu $.
The measure $\mu $ on $G$ is $\sigma $-finite,
since $0<\mu (W)<\infty $ and $G$ is with a countable base and a locally finite
covering as in \S 2.8 and
\S 2.9. A strong continuity of the regular representation
$T: G'\to U(L^2(G,\mu ,{\bf C}))$ follows from the continuity
of the quasi-invariance factor $\rho _{\mu }(\psi ,x)$ by $(\psi ,x)
\in G'\times G$ and the embedding $T_eG'\hookrightarrow T_eG$ of trace class,
where $T^{\mu }:=T$,
$T(z)f(g):=\rho _{\mu }^{1/2}(z,g)f(z^{-1}g)$, $z\in G'$, $g\in G$,
$f\in L^2(G,\mu ,{\bf C})$.
Let a $\nu $-measurable function $f:
H\to \bf C$ be such that $\nu ( \{ x\in H:  
f(x+y)\ne f(x) \} =0$ for each $y\in X_0$ with 
$f\in L^1(H,\nu ,{\bf C})$.  Let
also $P_k :  l_2 \to L(k)$ be projectors 
such that $P_k(x)=x_k$ for each
$x=(\sum_{j\in {\bf N}}x^je_j )$, where $x_k:=\sum_{j=1}^k x^je_j$,
$x_k\in L(k)$, $L(k):=sp_{\bf R}(e_1,...,e_k)$, 
$sp_{\bf R}(e_j: j\in {\bf N}):=\{ y:  y\in l_2; 
y=\sum_{j=1}^nx^je_j; x^j\in {\bf R};
n\in {\bf N} \} $.  Since the dense subspace $X$ in $H$
is isomorphic with $l_2$, then each finite-dimensional subspace $L(k)$
is complemented in $H$ \cite{nari}.
From the proof of
Proposition II.3.1 \cite{dal} in view of the Fubini Theorem there exists a
sequence of cylindrical functions $f_k(x)=f_k(x^k)=\int_{H\ominus L(k)}
f(P_kx+(I-P_k)y)\nu _{I-P_k}(dy)$ which converges to $f$ in 
$L^1(H,\nu ,{\bf
C})$, where $\nu =\nu _{L(k)}\otimes \nu _{I-P_k}$, $\nu _{I-P_k}$ is the
measure on $H\ominus L(k)$.  
Each cylindrical function $f_k$ is $\nu $-almost
everywhere constant on $H$, 
since $L(k)\subset X_o$ for each $k\in \bf N$,
consequently, $f$ is $\nu $-almost everywhere constant on $H$.  
From the construction of $G'$ and $\mu $ with the help of 
the local diffeomorphism $A$ and $\nu $ it follows
that, if a function $f\in L^1(G,\mu ,{\bf C})$ 
satisfies the following condition
$f^h(g)=f(g)$ $(mod $ $\mu )$ by $g\in G$ for each $h\in G'$, 
then $f(x)=const $
$( mod $ $\mu )$, where $f^h(g):=f(hg)$, $g\in G$.  
\par Let $f(g)=Ch_U(g)$ be
the characteristic function of a subset $U$, $U\subset G$, $U\in Af(G,\mu )$,
then $f(hg)=1 $ $\Leftrightarrow g\in h^{-1}U$.  If $f^h(g)=f(g)$ is true by
$g\in G$ $\mu $-almost everywhere, then $\mu (\{ g\in G:  f^h(g)\ne f(g) \}
)=0$, that is $\mu ( (h^{-1}U)\bigtriangleup U)=0$, consequently, 
the measure $\mu $ 
satisfies the condition $(P)$ from \S VIII.19.5 \cite{fell}, where
$A\bigtriangleup B:=(A\setminus B)\cup (B\setminus A)$ 
for each $A, B\subset G$.
For each subset $E\subset g_jW_j$ with $g_j\in G'$ and $W_j\subset W$ 
from \S 2.9 the outer measure is bounded,
$\mu ^*(E)\le 1$, since $\mu
(W)=1$ and $\mu $ is non-negative \cite{boui}, consequently, 
there exists $F\in
Bf(G)$ such that $F\supset E$ and $\mu (F)=\mu ^*(E)$.  
This $F$ may be interpreted as the least upper bound in $Bf(G)$ relative to
the latter equality.
In view of the
Proposition VIII.19.5 \cite{fell} the measure $\mu $ is ergodic, that is for
each $U\in Af(G,\mu )$ and $F\in Af(G,\mu )$ with $\mu (U)
\times \mu (F)\ne 0$
there exists $h\in G'$ such that $\mu ((h\circ E)\cap F)\ne 0$.  
\par From Theorem I.1.2 \cite{dal} 
it follows that $(G, Bf(G))$ is a Radon space,
since $G$ is separable and complete. Therefore, a class of compact subsets
approximates from below each measure $\mu ^f$, 
$\mu ^f(dg):=|f(g)|\mu (dg)$,
where $f\in L^2(G,\mu ,{\bf C})=:\bar H$.
Due to the Egorov Theorem 2.3.7 \cite{fed} for each $\epsilon >0$
and for each sequence
$f_n(g)$ converging to $f(g)$ for $\mu $-almost every $g\in G$,
when $n\to \infty $, there exists a compact subset $\sf K$
in $G$ such that $\mu (G\setminus {\sf K})<\epsilon $ and
$f_n(g)$ converges on $\sf K$ uniformly by $g\in \sf K$,
when $n\to \infty $.
In each Hilbert space $L^2({\bf R^n},\lambda ,{\bf R})$
the linear span of functions 
functions $f(x)=exp[(b,x)-(ax,x)]$ is dense, where $b$ and $x\in 
\bf R^n$, $a$ is a symmetric positive definite real $n\times n$ matrix,
$(*,*)$ is the standard scalar product in $\bf R^n$
and $\lambda $ is the Lebesgue measure on $\bf R^n$.
If a non-linear operator $U$ on $X$ satisfies conditions of Theorem 
26.1 \cite{sko}, then $\nu ^U(dx)/\nu (dx)=|det U'(U^{-1}(x))|
\rho _{\nu }(x-U^{-1}(x),x)$, where $\nu ^U(B):=\nu (U^{-1}B)$
for each $B\in Bf(X)$, $\rho _{\nu }(z,x)=exp\{ \sum_{l=1}^{\infty }
[2(z,e_l)(x,e_l)-(z,e_l)^2]/\lambda _l \}$ by Theorem 26.2 \cite{sko}, 
where $\lambda _l$ and $e_l$
are eigenvalues and eigenfunctions of the correlation operator $J$ on $X$
enumerated by $l\in \bf N$,
$z\in X_0$, $\rho _{\nu }(z,x):=\nu _z(dx)/\nu (dx)$, 
$\nu _z(B):=\nu (B-z)$ for each $B\in Bf(X)$. Since the Gaussian measure
$\nu $ induces with the help of subalgebras of cylinder subsets
in $Bf(H)$ and $Bf(X)$ the corresponding Gaussian measure
on $H$, which is also denoted by $\nu $,
then analogous formulas of quasi-invariance factor
are true for $\nu $ on $H$ \cite{dal}.
Hence in view of the Stone-Weierstrass Theorem A.8 \cite{fell}
an algebra ${\sf V}(Q)$ of finite pointwise products of 
functions from the following space 
$sp_{\bf C}\{ \psi (g):=(\rho (h,g))^{1/2}: h\in G' \}=:Q$ is dense in
$L^2(G,\mu ,{\bf C})$, since $\rho (e,g)=1$ for each $g\in G$
and $L_h: G\to G$ are diffeomorphisms of the manifold $G$, $L_h(g)=hg$.
\par For each $m\in \bf N$ there are $C^{\infty }$-curves
$\phi _j^b$, where $j=1,...,m$ and $b\in (-2,2):= \{ a: a\in {\bf R};
-2<a<2 \} $ is a parameter, such that $\phi _j^b|_{b=0}=e$
and $(\partial \phi _j^b/\partial b)|_{b=0}$ are linearly independent in
$T_eG'$ vectors and $\phi _j:=\phi _j^1$,
$\phi _j\in G'\cap W$, $j=1,...,m$, since $G'$ 
is the infinite-dimensional group, which is complete relative to its
own uniformity. 
Then the following condition
$det (\Psi (g))=0$ defines a submanifold $G_{\Psi }$
in $G$ of codimension over $\bf R$, 
\par $(i)$ $codim_{\bf R}G_{\Psi }\ge 1$, 
where $\Psi (g)$ is a matrix dependent from $g\in G$ with matrix elements
$\Psi _{l,j}(g):=D^{2l}_{\phi _j}(\rho (\phi _j,g))^{1/2}$.
If $f\in \bar H$ is such that 
$(f(g),(\rho (\phi ,g))^{1/2})_{\bar H}=0$
for each $\phi \in G'\cap W$, then differentials of these scalar products
by $\phi $ are zero. But ${\sf V}(Q)$ is dense in 
$\bar H$ and in view of condition $(i)$ this means that $f=0$,
since for each $m$ there are $\phi _j\in G'\cap W$ such that
$det \Psi (g)\ne 0$ $\mu $-almost everywhere on $G$, $g\in G$.
If $\| f\|_{ \bar H}>0$, then $\mu (supp(f))>0$,
consequently, $\mu ((G'supp(f))\cap W)=1$, since $G'U=G$ for each open
$U$ in $G$ and for each $\epsilon >0$ there exists an open $U$,
$U\supset supp(f)$, such that $\mu (U\setminus supp(f))<\epsilon $.
\par This means that the linear span over $\bf C$: $sp_{\bf C}
\{ Ch_{g_kW_k}\phi (g): \phi (g)=\rho ^{1/2}_{\mu }(h,g); h\in G' \} $
is dense in $L^2(g_kW_k,\mu ,{\bf C})$. 
Therefore, the following vector $f_0(g):=\sum_{j=0}^{\infty }2^{-j}Ch_{g_jW_j}
(g)$ is cyclic for $T^{\mu }$, since $\{ g_jW_j: j\in {\bf N_o} \} $
is a locally finite covering and ${\tilde \mu }(dg)=f_0(g)\mu (dg)$
is a finite measure with continuous $\rho _{\tilde \mu }$ such that
$f(g)\mapsto f^{1/2}_0(g)f(g)$ establishes isomorphism of $L^2(G,\tilde \mu ,
{\bf C})$ with $L^2(G,\mu ,{\bf C})$.
If $f_k\in L^{\infty }(g_kW_k,\mu ,{\bf
C})$ for each $k\in \bf N$ and $f_k|_{(g_kW_k\cap g_lW_l)}$ 
$=f_l|_{(g_kW_k\cap g_lW_l)}$ for each $g_kW_k\cap g_lW_l\ne \emptyset $
and $\sup_k\| f_k\|_{L^{\infty }(g_kW_k,\mu ,
{\bf C})}<\infty $, then there exists $f\in L^{\infty }(G,
\mu ,{\bf C})$ such that $f|_{g_kW_k}=f_k$ for each $k\in \bf N$, 
where $Ch_W(g)$ is the characteristic function of $W$, that is,
$Ch_W(g)=1$ for each $g\in W$ and $Ch_W(g)=0$ for each $g\in G\setminus W$. 
From the construction of $\mu $ it follows that for each $f_{1,j}$ and 
$f_{2,j}\in \bar H$, $j=1,...,n$, $n\in \bf N$ and each $\epsilon >0$ there 
exists $h\in G'$ such that $|(T_hf_{1,j},f_{2,j})_{\bar H}|\le 
\epsilon |(f_{1,j},f_{2,j})_{\bar H}|$,
when $|(f_{1,j},f_{2,j})_{\bar H}|>0$,
since $G$ is the Radon space by Theorem I.1.2 \cite{dal}
and $G$ is not locally compact. 
This means that there is not any finite-dimensional 
$G'$-invariant subspace $H'$ in $\bar H$ such that
$T_hH'\subset H'$ for each $h\in G'$ and $H'\ne \{ 0 \}$.
Hence if there is a $G'$-invariant closed subspace $H'$
in $\bar H$ it is isomorphic with the subspace
$L^2(V,\mu ,{\bf C})$, where $V\in Bf(G)$. 
\par Let ${\sf A}_G$ denotes a $*$-subalgebra of $L(\bar H, \bar H)=L(\bar H)$
generated by the family of unitary operators 
$\{ T_h: h\in G' \} $. In view of the von Neumann
double commuter Theorem (see \S VI.24.2 \cite{fell})
${{\sf A}_G}"$ coincides with the weak and strong operator closures of
${\sf A}_G$ in $L({\bar H}, {\bar H})$, where ${{\sf A}_G}'$
denotes the commuting algebra of ${\sf A}_G$ and ${{\sf A}_G}"=
({{\sf A}_G}')'$. Suppose that $\lambda $ is a probability 
Radon measure on $G'$ such that $\lambda $ has not any atoms and
$supp (\lambda )=G'$.
Let $a(x)\in L^{\infty }(G,\mu ,{\bf C})$, $f$ and $g\in \bar H$, 
$\beta (h)\in L^2(G',\lambda ,{\bf C})$. Since $L^2(G',\lambda ,{\bf C})$ 
is infinite-dimensional, then for each finite family of 
$a\in \{ a_1,...,a_m \} \subset L^{\infty }(G,\mu ,{\bf C})$,
$f\in \{ f_1,...,f_m \} \subset \bar H$ there exists
$\beta (h)\in L^2(G',\lambda ,{\bf C})$, $h\in G'$, such that
$\beta $ is orthogonal to $\int_G{\bar f}_s(g)
[f_j(h^{-1}g)
(\rho (h,g))^{1/2}-f_j(g)]\mu (dg)$ for each $s,j=1,...,m$. Hence
each operator of multiplication on $a_j(g)$
belongs to ${{\sf A}_G}"$, since there exists $\beta (h)$
such that $(f_s,a_jf_l)=$ 
$\int_G\int_{G'}{\bar f}_s(g)\beta (h)(\rho (h,g))^{1/2}f_l(h^{-1}g)
\lambda (dh) \mu (dg)$ $=\int_G\int_{G'} {\bar f}_s(g)
\beta (h)(T_hf_l(g))\lambda (dh)\mu (dg)$,
$\int_G{\bar f}_s(g)a_j(g)f_l(g)\mu (dg)$ 
$=\int_G \int_{G'}
{\bar f}_s(g)$ $\beta (h)$ $f_l(g)\lambda (dh)\mu (dg)=$
$(f_s,a_jf_l)$. Hence ${{\sf A}_G}"$ contains 
subalgebra of all operators of multiplication on functions from
$L^{\infty }(G,\mu ,{\bf C})$.
\par Let us remind the following. A Banach bundle $\sf B$ over 
a Hausdorff space $G'$ is a bundle $<B,\pi >$ over $G'$, together 
with operations and norms making each fiber $B_h$ ($h\in G'$)
into a Banach space such that conditions $BB(i-iv)$ are satisfied: 
$BB(i)$ $x\to \| x\| $ is 
continuous on $B$ to $\bf R$; $BB(ii)$ the operation $+$ is 
continuous as a function on $\{ (x,y)\in B\times B: 
\pi (x)=\pi (y) \} $ to $B$; $BB(iii)$ for each $\lambda \in \bf C$,
the map $x\to \lambda x$ is continuous on $B$ to $B$; $BB(iv)$ if $h\in G'$
and $\{ x_i\} $ is any net  of elements of $B$ such that $\| x_i\| 
\to 0$ and 
$\pi (x_i)\to h$ in $G'$, then $x_i\to 0_h$ in $B$, 
where $\pi : B\to G'$ is a bundle projection, 
$B_h:=\pi ^{-1}(h)$ is the fiber over $h$ (see \S II.13.4
\cite{fell}). If $G'$ is a Hausdorff topological group, then a Banach 
algebraic bundle over $G'$ is a Banach bundle ${\sf B}=<B,\pi >$ over $G'$
together with a binary operation $\bullet $ on $B$ satisfying
conditions $AB(i-v)$: 
$AB(i)$ $\pi (b\bullet c)=\pi (b)\pi (c)$ for $b$ and $c\in B$;
$AB(ii)$ for each $x$ and $y\in G'$ the product $\bullet $
is bilinear on $B_x\times B_y$ to $B_{xy}$;
$AB(iii)$ the product $\bullet $ on $B$ is associative;
$AB(iv)$ $\| b\bullet c\| \le \| b\| \times \| c\| $
($b, c\in B$); $AB(v)$ the map $\bullet $ is continuous on $B\times B$ 
to $B$ (see \S VIII.2.2 \cite{fell}). With $G'$ and a Banach algebra $\sf A$ 
the trivial Banach bundle ${\sf B}={\sf A}\times G'$ is associative, in 
particular let ${\sf A}=\bf C$ (see \S VIII.2.7 \cite{fell}).
\par The regular representation $T$ of $G'$ gives rise to a canonical regular
$\bar H$-projection-valued measure $\bar P$:
$\bar P(W)f=Ch_Wf$, where $f\in \bar H$, $W\in Bf(G)$, $Ch_W$ 
is the characteristic function of $W$. Therefore, $T_h\bar P(W)=\bar P
(h\circ W)T_h$ for each $h\in G'$ and $W\in Bf(G)$, since
$\rho (h,h^{-1}\circ g)\rho (h,g)=1=\rho (e,g)$ for each $(h,g)
\in G'\times G$, 
$Ch_W(h^{-1}\circ g)=Ch_{h\circ W}(g)$ and $T_h(\bar P(W)f(g))=\rho (
h,g)^{1/2}\bar P(h\circ W)f(h^{-1}\circ g)$. Thus $<T,\bar P>$ is 
a system of imprimitivity for $G'$ over $G$, which is denoted 
${\sf T}^{\mu }$. This means that conditions
$SI(i-iii)$ are satisfied: $SI(i)$ $T$ is a unitary representation
of $G'$; $SI(ii)$ $\bar P$ is a regular 
$\bar H$-projection-valued Borel measure on $G$ and 
$SI(iii)$ $T_h\bar P(W)=\bar P(h\circ W)T_h$ for all $h\in G'$ 
and $W\in Bf(G)$. 
\par For each $F\in L^{\infty }(G,\mu ,{\bf C})$ let $\bar \alpha _F$
be the operator in $L({\bar H}, {\bar H})=L(\bar H)$ consisting
of multiplication by $F$: $\bar \alpha _F(f)=Ff$, $f\in \bar H$. 
The map $F\to \bar \alpha _F$ is  an isometric $*$-isomorphism
of $L^{\infty }(G,\mu ,{\bf C})$ into $L({\bar H}, {\bar H})$
(see \S VIII.19.2\cite{fell}). Therefore, Propositions 
VIII.19.2,5\cite{fell}
(using the approach of this particular case given above) are applicable
in our situation.
\par If $\bar p$ is a projection onto a closed ${\sf T}^{\mu }$-stable
subspace of $\bar H$, then $\bar p$ commutes with all
$\bar P(W)$. Hence $\bar p$ commutes with multiplication by all
$F\in L^{\infty }(G,\mu ,{\bf C})$, so by VIII.19.2 \cite{fell}
$\bar p=\bar P(V)$, where $V\in Bf(G)$. Also $\bar p$ commutes with all
$T_h$, $h\in G'$, consequently, $(h\circ V)\setminus V$ and 
$(h^{-1}\circ V)\setminus V$ are $\mu $-null for each $h\in G'$, 
hence $\mu ((h\circ V)\bigtriangleup V)=0$ for all $h\in G'$. In view 
of ergodicity of $\mu $ and proposition VIII.19.5 \cite{fell}
either $\mu (V)=0$ or $\mu (G\setminus V)=0$, hence
either $\bar p=0$ or $\bar p=I$, where $I$ is the unit operator.
Hence $T$ is the irreducible unitary representation.
\par Almost analogous proof was done in the case of loop groups
with the corresponding quasi-invaraint measures 
and with the use of the spectral theorem  for the family of commuting unitary
operators, since the loop group is Abelian
in \cite{luihlgr,luihlgna}.
In the non-Archimedean case $G'$ has the analytic atlas
$At(G')=\{ (U_j,\psi _j): j \in {\bf N} \}$ with disjoint clopen charts, 
hence curves $\phi _j^b$ can be chosen locally
analytic with a restriction on the corresponding neighbourhood $U_1$ of $e$
being analytic, where $b\in \bf L$.
Substitution of differentiation on pseudodifferentiation along $\phi _j^b$
by parameter $b\in B({\bf L},0,1)$ produces by formula $det (\Psi (g))=0$
an analytic submanifold $G_{\Psi }$ in $G$ 
with $codim_{\bf L}G_{\Psi }\ge 1$, since $G$ is the analytic manifold.
\par {\bf 3.9. Theorem.} {\it Let $P_m $ 
be the ergodic Poisson measure on 
${\tilde \Gamma }_X$ as in \S 2.4, 2.9 
and $q$ be an irreducible representation of the
symmetric group $\Sigma _n$ ($q=I$ for $X=G$ and may be non-trivial
for $X=M$ finite-dimensional over the coresponding field $\bf L$
and a group of diffeomorphisms $G'$ of $M$). Then the representation
$U_m \otimes V^q_m$ from \S 3.1 is irreducible.}
\par {\bf Proof.} The case of real finite-dimensional $M$ was proved
in \cite{vgg}. The case of non-Archimedean $M$ with $dim_{\bf L}M<\infty $
follows from \S 3.5, since ${\tilde \Gamma }_X=\Gamma _X$ in this case.
Indeed, $U_m\otimes V^q_m$ is equivalent with $U^q_{n\circ m}$
and $f\phi \in L^2(\Gamma _X,P_m,{\bf C})$, if $f\in L^2(\Gamma _X,P_m,{\bf C})$
and $\phi \in L^{\infty }(\Gamma _X,P_m,{\bf C})$. Then each subspace
$\sf L$ in $L^2(\Gamma _X,P_m,{\bf C})\otimes H^q$ 
invariant relative to $G'=Diff^t(X)$ is also invariant relative to 
multiplications on functions $\phi \in L^{\infty }(\Gamma _X,P_m,{\bf C})$,
since ${\sf L}=\bigoplus_{r,i}{\sf L}^i_{l,r}$, where
${\sf L}^i_{l,r}:={\sf L}\cap (L^2(B^r_{K_l}\times \Gamma _{X\setminus K_l},
\mu _r,{\bf C})\otimes W^i_r\otimes C^i_r)$ are subspaces invariant relative to
$Diff^t(K_l)$, $\mu =P_m$ and $\mu _r$ is the corresponding measure on
$B^r_{K_l}\times \Gamma _{X\setminus K_L}$. 
In view of lemma 2.14 the measure $\mu _n$ is equivalent with
${\mu '}_n\times {\mu "}_n$ and further as at the end of \S 3 \cite{vgg}.
\par The remaining cases are proved analogously to the proof of theorem
3.8 (and see \cite{lutmf,luihlgr,luihlgna})
applied to the pair $(G',{\tilde \Gamma }_X)$
instead of $(G',G)$, since ${\tilde \Gamma }_X$ is $C^{\infty }$-manifold
and from infinite differentiability or pseudodifferentiability 
of $m$ it follows,
that $P_m$ is also infinite differentiable or pseudodifferentiable respectively,
morever, $P_m$ is the ergodic measure due to theorem 2.9.
In the case of $X=M$ the measures on $X$ are chosen to be such that
$sp_{\bf C}\{ \rho _m^{1/2}(z,x)=\phi (x): z\in G' \} $ is dense
in $L^2(X,m,{\bf C})$ in accordance with \S 2.9 and the cited papers there, 
for example, Gaussian measures or product measures of special type
on $T_xM$ induce the demanded measures on $M$,
where $x\in M$.
\par It remains only to establish that the density $\rho _{P_m}$
has the demanded properties. For this 
it is necessary to use the fact that operators $L_h$ on $X$ (either $X=M$
or $X=G$) are infinitely strongly differentiable by $h\in G'$
and there exists a dense subset $G"$ in $G'$ such that $(L_h)^{(n)}\ne 0$
for each $n\in \bf N$ and each $h\in G"$. Therefore, $(AL_hA^{-1})^{(n)}\ne 0$
for each $h\in G"$, where $A: U\to V_H$ is a local diffeomorphism, where
$U$ and $V_H$ are open subsets in $X$ and the corresponding Banach space $H$
respectively as in \S 2.9, \S 3.8 and the cited above papers.
In the Hilbert space $L^2({\bf R^{mn}},\lambda ,{\bf C})$
is dense the following linear span
$sp_{\bf C} \{ exp[\sum_{l=1}^k(a^l_1,x^l)-(a^l_2x^l,x^l)]=:\phi (y)|
a^l_1\in {\bf R^m},$ $a^l_2\in {\bf R^m}$, 
$\sum_{l=1}^ka^l_{2,j}>0$ for each $j=1,...,m;$  $a^l_{2,j}\ge 0$
for each $j=1,...,m,$ $l=1,...,k$; and if $a^l_{1,j}\ne 0$, then
there exists $l'$ such that $2l'>l$ and $a^{l'}_{2,j}>0$,
$a^l_i:=(a^l_{1,i},...,a^l_{i,m})$, 
$a^l_{i,j}\in \bf R$, $x_j^l:=S_l(y_{j,1},...,
y_{j,n})$, $l=1,...,k$, $j=1,...,m \} $ , where $k=k(m,n)\in \bf N$
is chosen 
such that $z\mapsto \{ S_l(z): l=1,..,k \} $ is a bijection of
$\bf R^n$, $z=(z_1,...,z_n)\in \bf R^n$, 
$S_l(z):=\sum_{i=1}^n(z_i)^l$ is a power sum
of degree $l$, $(*,*)$ is a scalar product in $\bf R^m$, $\lambda _{mn}$
is a Lebsgue measure on $\bf R^{mn}$.
For the local field $\bf L$
in the Hilbert space $L^2({\bf L^{mn}},\nu ,{\bf C})$
is dense the following linear span
$sp_{\bf C} \{ exp[\sum_{l=1}^k -|(a^l,(b^l+x^l)|^2]=:\phi (y)|
a^l\in {\bf L^m},$ $b^l\in {\bf L^m}$, 
$\sum_{l=1}^k|a^l_j|>0$ for each $j=1,...,m$, 
$a^l:=(a^l_1,...,a^l_m)$, $a^l_j\in \bf L$, $x_j^l:=S_l(y_{j,1},...,
y_{j,n})$, $l=1,...,k$, $j=1,...,m \} $ , where $k=k(m,n)\in \bf N$
is chosen  
such that $z\mapsto \{ S_l(z): l=1,..,k \} $ is a bijection of
$\bf L^n$, $z=(z_1,...,z_n)\in \bf L^n$, 
$S_l(z):=\sum_{i=1}^n(z_i)^l$ is a power sum
of degree $l$, $(z,q):=\sum_{i=1}^nz_iq_i$, $q\in \bf L^n$, $\nu _{mn}$
is the Haar measure on $\bf L^{mn}$.
\par Using charts in $B^n_X$ we get projections $L^2(B^n_X,m^n,{\bf C})$
into $L^2({\bf R^{mn}},\lambda _{mn},{\bf C})$ in the real case
and into $L^2({\bf L^{mn}},\nu _{mn},{\bf C})$ in the non-Archimedean case.
Then we use Taylor expansion up to $o(d^{k'+1}_{G'}(\psi ,h))$
of $L_h$ in a suitable neighbourhoods $hU'$ in $G'$
of elements $h\in G"$ with $U'$ open in $G'$ with $e\in U'\subset W'$
and with $k'=2k$ in the real case and $k'=k$ in the non-Archimedean case,
where $d_{G'}$ is the metric in $G'$ in its own uniformity, $\psi
\in hU'$. For a manifold $C^t(M,N)$ of mappings $f: M\to \bf N$
of class of smoothness $C^t$ with $t\ge 1$
for $C^{\infty }$-manifolds $M$ and $N$
in the real case and analytic manifolds in the non-Archimedean case
the tangent manifold $TC^t(M,N)$ is isomorphic with $C^t(M,TN)$
and for the $n$-th order we get $T^nC^t(M,N)=C^t(M,T^nN)$ (see also
\cite{ebi,kling,luihlgna}).
Then $T^nDiff^t(M)$ is a submanifold in $C^t(M,T^nM)$.
Let $C^t(M,m_0;N,n_0)$ be a family of mappings $f\in C^t(M,N)$ preserving 
marked points $f(m_0)=n_0$, $m_0\in M$ and $n_0\in N$, where
in the real case $M=S^m$ is the $m$-dimensional real sphere and 
$dim_{\bf R}N>m$. Analogously for others classes of smoothness $\omega $
considered for construction of loop groups $L(M,m_0;N,n_0)_{\omega }$, 
elements of which
are closures of orbits $cl \{ f(\psi (x)): \psi \in G(M), \psi (m_0)=m_0 \} $,
where $G(M)$ denotes the group of diffeomorphisms of $M$ of the corresponding
class of smoothness and with certain additional construction in the
non-Archimedean case \cite{luumnad,luihlgr,luihlgna}. 
Hence the manifold $T^nL(M,m_0;N,n_0)_{\omega }$ is isomorphic with 
the following manifold $T^nL(M,m_0;T^nN,
(n_0,{\bar 0}_n))_{\omega }\otimes T^n_{(n_0,{\bar 0}_n)}N$, where
${\bar 0}_n\in T^nN$ is the zero section for each $n\in \bf N$.
Therefore, it is possible to vary
values of differentials $D^jf$ for $j=0,...,n$
in the notation $T^nf:=(f,Df,D^2f,...,D^nf)$ with $D^0f:=f$
for elements $f\in G"$ both in the case of the 
diffeomorphism group 
and the loop group up to the corrections $o(d^{n+1}_{G'}(f,\psi ))$.
Then $D^n_hL_h(g)$ can be expressed through
$D^jh$ and $D^jg$ with $1\le j\le n$, where $h, g\in G$,
hence it is possible to vary coefficients $a^l_1$, $a^l_2$ in the real case
and $a^l$, $b^l$ in the non-Archimedean case.
\par Take for example, the Gaussian measure on $X$ in the real case
induced from the Gaussian measure on the corresponding Banach space
and given with the help of non-degenerate symmetric positive definite
operator of trace class. In the non-Archimedean case each Banach space over
a local field $\bf L$ is isomorphic with $c_0(\alpha ,{\bf L})$, 
where $\alpha $ is an
ordinal and elements of $c_0(\alpha ,{\bf L})$ have the form $x=(x_j: j\in
\alpha , x_j\in {\bf L})$ such that $\| x\| :=\sup_j|x_j|<\infty $
and for each $\epsilon >0$ a set $\{ j: |x_j|>\epsilon \} $ is finite
\cite{roo}.
For each separable manifold $M$ we have $card(\alpha )\le \aleph _0$.
In the latter case take, for example, the following non-Archimedean analog
$\eta $ of the Gaussian measure: each projection $\eta _j$
of $\eta $ on ${\bf L}e_j$ has a density $\eta _j(dx)=F_j exp(-|x|^2s_j)
v(dx)$, where $\sum_js_j^{-1}<\infty $, $e_j:=(0,...,0,1,0...)$
with $1$ on the $j$-th place,
$v$ is the Haar measure on $\bf L$
with $v(B({\bf L},0,1))=1$ and constants $F_j>0$ are chosen such that
$\eta _j({\bf L})=1$ (see also \S 2.9).
\par Let $\psi _j^b$ be $C^{\infty }$-curves 
in $G'$ such that $(\partial \psi _j^b/\partial b)|_{b=0}$ 
are linearly independent vectors
in $T_{g_k}G'$ and $\psi _j^b|_{b=0}=g_k$ and 
$R(\psi _j)\cap R(\psi _l)\cap (g_kW_k)
=\{ g_k \}$ for each $j\ne l$, $b\in \bf L$, 
where either $X=M$ or $X=G$, $j=1,...,n$, $n\in \bf N$,
$R(\psi )$ denotes the range of 
$\psi $, that is, $R(\psi ):=\{ \psi ^b: b\in {\bf L} \} $, $g_k\in G'$,
$S_k$ is open in ${\tilde \Gamma }_X$, $\zeta _k: S_k\to Q_k$
are local diffeomorphisms of open subsets $S_k$ in ${\tilde \Gamma }_X$
and $\Gamma _X$ respectively, $\gamma \in g_kS_k$ (see also \S 2.4).
\par There are embeddings $L^2(B^n_K,m_n,{\bf C})\hookrightarrow
L^2(B^n_X,m_n,{\bf C})\hookrightarrow L^2(\Gamma _X,P_m,{\bf C})
\hookrightarrow L^2({\tilde \Gamma }_X,P_m,{\bf C})$, where $m_n$ denotes
the restriction of $P_m$ on $B^n_X$. For each $x\in X$ there exists
$K\in \{ K_l: l\in {\bf N} \} $ such that $x\in Int(K)$. Then
$\Gamma _K$ is the disjoint union of $\{ B^n_K: n\in {\bf N_o} \} $.
On the other hand, $P_{K,m}|_{Bf(K^n)}=m_{K,n}$ in accordance with \S 3.1,
where $m_{K,n}$ is equivalent with $m^n_K$, where $m_K$ denotes the restriction
of $m$ on $Bf(K)$ and $m^n_K$ is the product of $n$ copies of $m_K$.
Then the condition
$det (\Psi (\gamma ))=0$ defines a submanifold ${\tilde \Gamma }_{X,\Psi }$
in ${\tilde \Gamma }_X$ of codimension over $\bf L$, 
\par $(i)$ $codim_{\bf L}{\tilde \Gamma }_{X,\Psi }\ge 1$, 
where $\Psi (\gamma )$ 
is a matrix dependent from $\gamma \in {\tilde \Gamma }_X$ with indices
of rows and columns $j$ and $l=1,...,n$
for $n\in \bf N$ with matrix elements
$\Psi _{l,j}(\gamma ):=D^{2l}_{\phi _j}(\rho (\phi _j,g))^{1/2}$
in the real case and with the corresponding pseudodifferentials 
by parameters $b_j\in B({\bf L},0,1)$ and for $\phi _j^{b_j}$
in the non-Archimedean case instead of differentials (see also
\cite{lutmf,lupr210}). 
In the non-Archimedean case ${\tilde \Gamma }_X$ has the analytic atlas
$At({\tilde \Gamma }_X)=\{ (V_j,\omega _j): j \in {\bf N} \}$ 
with disjoint clopen charts, also $G'$ has disjoint clopen charts
and the analytic atlas $At(G')=\{ (U_j,\phi _j): j \} $,
hence curves $\psi _j^b$ can be chosen locally
analytic with a restriction on the corresponding neighbourhood $U_1$ of $e$
in $G'$ being analytic, where $b\in \bf L$.
Substitution of differentiation on pseudodifferentiation along $\phi _j^b$
by parameter $b\in B({\bf L},0,1)$ produces by formula $det (\Psi (\gamma ))=0$
an analytic submanifold ${\tilde \Gamma }_{X,\Psi }$ 
in ${\tilde \Gamma }_X$ with $codim_{\bf L}{\tilde \Gamma }_{X,\Psi }\ge 1$.
Since for equivalent measures such regular
representations are equivalent, we can consider infinitely
differentiable or pseudodifferentiable measures in the real and non-Archimedean
cases respectively. There is the following equality 
$\lim_{(m^n(B)\to 0, \infty >m^n(B)>0)}m^n(B)exp(-m^n(B))/
(({m^n})^{\psi }(B)exp(-({m^n})^{\psi }(B)))$ $=:m^n(dx)exp(-m^n(dx))/
(({m^n})^{\psi }(dx)exp(-({m^n})^{\psi }(dx)))$ $=\rho _{m^n}(\psi ,x)$,
since $\rho _{m^n}(\psi ,x)$ is continuous on $G'\times {\tilde X}^n$ and
$\lim_{(m^n(B)\to 0, \infty >m^n(B)>0, x\in B)}exp(-{m^n}(B)[1-\int_B
\rho _{m^n}(\psi ,y)m^n(dy)/m^n(B)])=1$ for balls $B$ in ${\tilde X}^n$
such that $x\in Int(B)$, where $({m^n})^{\psi }(E):=m^n(\psi ^{-1}E)$
for each Borel subset $E\in Bf({\tilde X}^n)$, $\psi \in G'$,
$x\in {\tilde X}^n$.
\par In the case $X=G$ each space $L^2(B^r_X,m_r,{\bf C})$
has the embedding into $L^2({\tilde \Gamma }_X,P_m,{\bf C})$,
where $m_r:=P_m|_{B^r_X}$.
It was supposed above that the quasi-invariance factor $\rho _m(\psi ,x)$
of the quasi-invariant measure $m$ on $Bf(X)$ relative to $G'$
is continuous on $G'\times X$, consequently, $\rho _{m^r}(\psi ,\eta )$
and $\rho _{m_r}(\psi ,\gamma ^r)$ and 
$\rho _{P_m}(\psi ,\gamma )$ are continuous on $G'\times X^r$
and $G'\times B^r_X$ and $G'\times \Gamma _X$ respectively, where
$\psi \in G'$, $\eta \in X^r$, $\gamma ^r\in B^r_X$, $\gamma \in \Gamma _X$.
Hence due to the definiton of $P_m$ there is the equality:
$lim_{r\ge n, r\to \infty }\rho _{m_r}(\psi ,\gamma ^n)=
\rho _{P_m}(\psi ,\gamma ^n)$ for each $\psi \in G'$,
$\gamma ^n\in B^n_X$.
If $f\in \bar H:=L^2({\tilde \Gamma }_X,P_m,{\bf C})$ is such that 
$(f(g),(\rho _{P_m}(\phi ,g))^{1/2})_{\bar H}=0$
for each $\phi \in G'\cap W$, then differentials of these scalar products
by $\phi $ are zero. In view of the above embeddings and formula
2.4$(i)$ and in view of condition $(i)$ this means that $f=0$,
since for each $n\in \bf N$ there are $\phi _j\in G'\cap W$ such that
$det \Psi (\gamma )\ne 0$ $P_m$-almost everywhere on ${\tilde \Gamma }_X$, 
$\gamma \in {\tilde \Gamma }_X$.
If $\| f\|_{ \bar H}>0$, then $P_m(supp(f))>0$,
consequently, $P_m(G'supp(f))=1$, 
since $G'U={\tilde \Gamma }_X$ for each open
$U$ in ${\tilde \Gamma }_X$ and for each $\epsilon >0$ there exists an open $U$,
$U\supset supp(f)$, such that $P_m(U\setminus supp(f))<\epsilon $.
\par This means that the linear span over $\bf C$: $sp_{\bf C}
\{ Ch_{g_kS_k}\phi (g): \phi (g)=\rho ^{1/2}_{P_m}(h,g); h\in G' \} $
is dense in $L^2(g_kS_k,P_m,{\bf C})$, since ${U'}_lK_l\subset Int(K_{l+1})$
for each $l\in \bf N$ (see \S 2.9). 
Therefore, $sp_{\bf C} \{ \phi (g): \phi (g)=\rho ^{1/2}_{P_m}(h,g); 
h\in G' \} $ is dense in $L^2({\tilde \Gamma }_X,P_m,{\bf C})$
and a vector $f_0$ is cyclic for $U_m$, where
$f_0(\gamma )=1$ for each $\gamma \in {\tilde \Gamma }_X$.
Then ${{\sf A}_G}"$ contains 
subalgebra of all operators of multiplication on functions from
$L^{\infty }({\tilde \Gamma }_X,P_m,{\bf C})$ and the remainder of the proof 
of theorem 3.9 is quite analogous with the proof of theorem 3.8
(certainly ${{\sf A}_G}"\ne L^{\infty }({\tilde \Gamma }_X,P_m,{\bf C})I$
for $G'=Diff^t(M)$,
since the regular representation $U_m(h)$
of $G'$ contains a family of
cardinality ${\sf c}:=card({\bf R})$ of non-commuting operators
from the set $ \{ U_m(h): h\in G' \} $).
\par {\bf 3.10. Theorem.} {\it $(\alpha ).$ 
If there exists a bounded operator
$T: L^2({\tilde \Gamma }_X,P_m,{\bf C})\otimes H^q\to 
L^2({\tilde \Gamma }_X,P_{m'},
{\bf C})\otimes H^{q'}$ ($H^q=\{ 0\} $ and $P_m$ is from theorem 3.9
for $X=G$ or infinite-dimensional
$X=M$ over the corresponding field $\bf L$ such that $L^2\otimes \{ 0\}:=L^2$)
satisfying conditions $(a,b)$:\\
$(a)$ $TU^q_m(\psi )=U^{q'}_{m'}T$ for all $\psi \in G'$,\\
$(b)$ there exists $\phi \in H^q$ such that $T(1\otimes \phi )\ne 0$,\\
then $P_m$ and $P_{m'}$ are equivalent.
\par $(\beta ).$ If there exists a bounded operator $V: L^2(G,\mu ,{\bf C})
\to L^2(G,\mu ',{\bf C})$ such that $VT^{\mu }(\psi )=T^{\mu '}
(\psi )V$ for each
$\psi \in G'$, where $\mu $ is a quasi-invaraint measure on $G$ relative to $G'$
and $T^{\mu }$ is the associated regular representation of $G'$ from theorem
3.8, then $\mu $ and $\mu '$ are equivalent.}
\par {\bf The proof} is divided into several parts. 
At first the case $(\alpha )$
of $X=M$ finite-dimensional over the corresponding filed $\bf
L$ is considered in subparagraphs {\bf I-III}. The cases $(\alpha )$
of $X=G$
and infinite-dimensional $X=M$ over $\bf L$ and the cases $(\beta )$
are considered in \S 3.10.{\bf IV}.
\par {\bf I.} Suppose that 
$\| \phi \| =1$ and $T$ is a contraction
operator. Take $X_n:=K_n$, where $n\in \bf N$ and $K_n$ are the same
as in \S 2.9. In the case $X=M$ we 
put $Y=X_n$, $\mu =P_m$, $\mu '=P_{m'}$, $\mu _1$ and $\mu _2$
are equal to the image measure of $\mu $ in accordance with the maps:
$\gamma \mapsto (\gamma \cap Y)=:\gamma _1$, $\gamma \mapsto (\gamma \cap
Y^c)=:\gamma _2$. 
Apart from the case $X=M$, for $X=G$ we suppose that
$Y=X$, since $G'$ acts on $G$ transitively and $supp(L_{\psi }):=cl \{
g\in G: \psi g\ne g \} =G$ for each $\psi \ne e$, because $G'$ is a dense
subgroup of $G$ and from $hg=g$ it follows $h=e$, where $h, g \in G$.
In the case of $Diff^t(X)$ there exists a bounded operator $T_Y: L^2(\Gamma _Y,
\mu _1,{\bf C})\otimes H^q\to L^2(\Gamma _Y,\mu '_1,{\bf C})\otimes
H^{q'}$ such that
\par $(i)$ $T_YF(\gamma ,a')=\int_{\Gamma _{Y^c}}TF(\gamma _1,\gamma _2,
a')\mu '_2(d\gamma _2)$. Then $L^2(\Gamma _Y,\mu _1,{\bf C})$
is embeddable as a closed subspace into $L^2(\Gamma _X,\mu ,{\bf C})$
by the map $L^2(\Gamma _Y,\mu _1,{\bf C})\ni f(\gamma )\mapsto
\hat f(\gamma ):=f(\gamma \cap Y)\in L^2(\Gamma _X,\mu ,{\bf C})$.
Therefore, $T_YF$ depends on $(\gamma _1,a')$ and $T_YF(\gamma ,a'_{\sigma })=
q'(\sigma )^{-1}T_YF(\gamma ,a')$ for all $\sigma \in \Sigma _{n'}$, where
$a'_{\sigma }=(i_{\sigma (1)},...,i_{\sigma (n')})$ for each
$a'=(i_1,...,i_{n'})\in {\tilde {\bf N}}^{n'}$. 
Then for $Diff^t(X)$: $\| T_YF\|^2=$
$\sum_{a'\in {\tilde {\bf N}}^{n'}}$ $\int_{\Gamma _Y}
\| T_YF(\gamma _1,a')\|^2_{W'}$
$\mu '_1(d\gamma _1)$ $\le 
\sum_{a'\in {\tilde {\bf N}}^{n'}}$ $\int_{\Gamma _X}
\| T_YF(\gamma ,a')\|^2_{W'}$
$\mu '_1(d\gamma )$ $\le 
\int_{\Gamma _Y}\int_{\Gamma _Y^c}$ 
$\sum_{a'\in {\tilde {\bf N}}^{n'}}
\| TF(\gamma _1,\gamma _2,a')\|^2_{W'}$ 
$\mu '_1(d\gamma _1)\mu '_2(d\gamma _2)$ $=\| TF\|^2$ $\le \| F\| ^2$
in the case $X=M$, consequently, $\| T_Y\| \le 1$
and $T_Y$ is a contraction too. When $\psi \in Diff^t(Y)$, 
then $\sigma (\psi ,
\gamma )$ is independent of $\gamma _2$. Hence $T_YU^q_m(\psi )=
U^{q'}_{m'}(\psi )T_Y$ for each $\psi \in Diff^t(Y)$.
\par There exists the decomposition of $\Gamma _X$ into disjoint union
of subsets $B^r_{X_k}\times \Gamma _{X\setminus X_k}$ for $r=0,1,2,...$,
where each such subset is invariant relative to $Diff^t(X_k)$, where $k$
is fixed and $B^r_{X_k}$ is the set of $r$-point subsets in $X_k$,
$B^0_{X_k}:=\{ \emptyset \} $ is the singleton. Therefore,
$L^2(\Gamma _X,\mu ,{\bf C})\otimes H=\bigoplus_{r=0}^{\infty }
(L^2(B^r_{X_k}\times \Gamma _{X\setminus X_k}, \mu _r,{\bf C})\otimes H$,
where $\mu _k$ is a restriction of $\mu $ on $B^r_{X_k}\times \Gamma
_{X\setminus X_k}$ (see also \S 2.13). Hence $L^2(B^r_{X_k}\times
\Gamma _{X\setminus X_k}, \mu _r,{\bf C})\otimes H=
\bigoplus_i(L^2(B^r_{X_k}\times \Gamma _{X\setminus X_k}, \mu _r,{\bf C})
\otimes W^i_r\otimes C^i_r),$ where $W^i_r$ are spaces in which irreducible
pairwise non-equivalent unitary representations $q^i_r$ of the symmetric group
$\Sigma _r$ act, $C^i_r$ denote spaces in which $\Sigma _r$ acts trivially.
Each term in the direct sum is invariant under $Diff^t(X_k)$ such that from
$\psi \in Diff^t(X_k)$ and $\gamma \in B^r_{X_k}\times \Gamma _{X\setminus X_k}$
it follows that $\sigma (\psi ,\gamma )\in \Sigma _r$.
\par In view of lemma 2.14 the measure $\mu _r$ is equivalent with
$\mu '_r\times \mu "_r$. Hence there exists the isomorphism $\tau _r:
L^2(B^r_{X_k}\times \Gamma _{X\setminus X_k},\mu _r,{\bf C})\to
L^2(B^r_{X_k},\mu '_r,{\bf C})\otimes L^2(\Gamma _{X\setminus X_k},\mu "_r,
{\bf C})$ given by the following formula: $\tau _rF:=
[\mu _r(d\gamma )/(\mu '_r(d\gamma ')\mu "_r(d\gamma "))]^{1/2}$ $F(\gamma )$.
Hence there exists isomorphism $\tau _r: L^2(B^r_{X_k}\times \Gamma _{X\setminus
X_k},\mu _r,{\bf C})\otimes W^i_r\otimes C^i_r\to L^2(B^r_{X_k},\mu '_r,{\bf C})
\otimes W^i_r\otimes L^2(\Gamma _{X\setminus X_k},\mu "_r,{\bf C})
\otimes C^i_r$. Therefore, for finite-dimensional
manifolds $M$ over $\bf R$ or the local field $\bf L$ considered here
for $\mu =P_m$
there is true the following lemma (for finite-dimensional real $M$ see
also lemma 3.2 \cite{vgg}).
\par {\bf 3.11. Lemma.} {\it Under the isomorphism $\tau _r$ the operator
$U(\psi ):=U^q_m (\psi )$ for each  $\psi \in Diff^t(X_k)$
transforms into $\tau _rU(\psi )\tau _r^{-1}=U^i_r(\psi )\otimes I$, where
$I$ is the unit operator in $L^2(\Gamma _{X\setminus X_k},\mu "_r,{\bf C})
\otimes C^i_r$ and $U^i_r$ is the operator in the space $L^2(B^r_{X_k},
\mu '_r,W^i_r)$ such that $(U^i_r(\psi )F)(\gamma ^r)=\rho _{\mu '_r}^{1/2}
(\psi ,\gamma ^r)$ $q^i_r(\sigma _r(\psi ,\gamma ^r))F(\psi ^{-1}\gamma ^r)$,
where $\gamma ^r\in B^r_{X_k}$.}
\par {\bf II.} A unitary representation $Q: \Sigma _{\infty }\to U(H^q)$
such that $Q(\sigma ): \phi (a) \mapsto \phi (\sigma ^{-1}a)$ restricted on
$\Sigma _r$ splits into the direct sum of subspaces: 
$H^q=\bigoplus_iW^i_r\otimes C^i_r$, that is, $Q(\sigma )\phi =\sum_i
\{ q^i_r\otimes id \} \phi _{r,i}$, $\phi =\sum_i\phi _{r,i}$,
where $\phi _{r,i}\in W^i_r\otimes C^i_r$, $q^i_r$ are the ireducible 
and pairwise distinct representations of $\Sigma _r$, $\Sigma _{\infty }:=
ind-\lim_r\Sigma _r$. Since $\Gamma
_Y=\bigcup_{r=0}^{\infty }B^r_Y$ is the disjoint union of $B^r_Y$, then
there are the following orthogonal decompositions:
$L^2(\Gamma _Y,\mu _1,{\bf C})=\bigoplus_{r=0}^{\infty }$
$L^2(B^r_Y,\mu _1,{\bf C})$ and $L^2(\Gamma _Y,\mu _1,{\bf C})\otimes H^q=
\bigoplus_{r,i}\phi _{\mu }(r,i)$, where $\phi _{\mu }(r,i):=
L^2(B^r_Y,\mu _1,{\bf C})\otimes W^i_r\otimes C^i_r$ are invariant subspaces
of the representation $U^q_m|_{Diff^t(Y)}$. Therefore, $U^q_m(\psi )=
U^{r,i}_{\mu }(\psi )\otimes id$ on $\phi _{\mu }(r,i)$, where
\par $(II.i)$ $U^{r,i}_{\mu }(\psi )(F\otimes w^i_r)=\rho _{\mu _1}^{1/2}
(\psi ,\gamma _1)F(\psi ^{-1}(\gamma _1))q^i_r(\sigma (\psi ,\gamma ))w^i_r$
for $F\in L^2(B^r_Y,\mu _1,{\bf C})$ and $w^i_r\in W^i_r$. The irreducible 
unitary representations $U^{r,i}_{\mu }$ and $U^{r',i'}_{\mu }$ are
equivalent if and only if $i=i'$ and $r=r'$. 
\par Hence there exists the unique integer $J_i$ such that 
either $T_Y\phi _{\mu }(r,i)=0$ or 
$T_Y\phi _{\mu }(r,i)
\subset \phi _{\mu '}(r,J_i)$ and the representations $q^i_r$ and 
${q'}^{J_i}_{r'}$ are equivalent, consequently, $J_i\ne J_k$ for 
each $i\ne k$. 
\par There exist intertwining operators $\omega _{r,i}: W^i_r\to {W'}^{J_i}_r$
of the representations $q^i_r$ and ${q'}^{J_i}_r$. We denote by $J_Y$
the unitary operator $J_Y: L^2(B^r_Y,\mu _1,{\bf C})\to
L^2(B^r_Y,{\mu '}_1,{\bf C})$  given by the following formula:
$J_YF(\gamma _1):=(\mu _1(d\gamma _1)/{\mu '}_1(d\gamma _1))^{1/2}F(\gamma _1)$.
Hence 
\par $(II.ii)$ $U^{r,J_i}_{\mu '}(\psi )T_{r,i}=T_{r,i}U^{r,i}_{\mu }(\psi )$
for each $\psi \in Diff^t(Y)$, where $T_{r,i}=J_Y\otimes \omega _{r,i}:
L^2(B^r_Y,\mu _1,{\bf C})\otimes W^i_r\to L^2(B^r_Y,{\mu '}_1,{\bf C})\otimes
{W'}^{J_i}_r$. 
\par {\bf III.} Using the general fact of the representation theory
of topological groups from \S III of the proof of theorem
3.1 in \cite{shim} we get for each $(r,i)$: either exists
a bounded operator $U_{r,i}: C^i_r\to {C'}^{J_i}_r$
such that $T_Y|_{\phi _{\mu }(r,i)}=T_{r,i}\otimes U_{r,i}$
or $T_Y\phi _{\mu }(r,i)=0$.
Hence for $Diff^t(Y)$ there is the following equality:
\par $(III.i)$ $T_Y(1\otimes \phi )(\gamma ,a')={\sum '}_{r,i}
T_{r,i}\otimes U_{r,i}(\chi _{B^r_Y}\otimes \phi _{r,i})(\gamma ,a')$
$=(mu _1(d\gamma _1)/{\mu '}_1(d\gamma _1)$ ${\sum '}_{r,i}\chi _{B^r_Y}
(\gamma _1)(\omega _{r,i}\otimes U_{r,i})(\phi _{r,i})(a')$, where
$\sum '$ is a sum for $(r,i)$ such that $T_Y\phi _{\mu }(r,i)\ne 0$,
$\phi =sum_i\phi _{r,i}$, $\phi _{r,i}\in W^i_r\otimes C^i_r$,
$\chi _A$ is the characteristic function of the subset $A$.
Then $\| {\sum '}_{r,i}\chi _{B^r_Y}(\gamma _1)(\omega _{r,i}\otimes U_{r,i})
(\phi _{r,i})(a') \| ^2_{W'}$ $\le \sum_r\chi _{B^r_Y}(\gamma _1)
\sum_i\| \phi _{r,i} \|^2=1$.
\par In the case of $Diff^t(X)$ if $P_m$ and $P_{m'}$ are mutually singular,
then $\lim_{k\to \infty }T_{X_k}(1\otimes \phi )(\gamma ,a')$ 
converges to $T(1\otimes \phi )
(\gamma ,a')$ for $P_{m'}$-a.e. $\gamma $ due to the martingale convergence
theorem, but $\lim_{k\to \infty }
(\mu _1(d\gamma _1)/{\mu '}_1(d\gamma _1))^{1/2}(\gamma _1){\sum '}_{r,i}
\chi _{B^r_Y}(\gamma _1)(\omega _{r,i}\otimes U_{r,i})(\phi _{r,i})(a')=0$
for $P_m$-a.e. $\gamma $, hence $T(1\otimes \phi )=0$, which contradicts
the assumption of this theorem. 
\par {\bf IV.} In view of theorem 3.9 representations $U_m$ are irreducible
for $X=G$ or infinite-dimensional $X=M$ over the field $\bf L$.
It was proved in \S 3.9 that 
\par $(IV.i)$ the weak closure of subalgebra
generated by the family $\{ U_m(h): h\in G' \} $
in the algebra of bounded linear operators $L(\bar H)$
contains all operators of multiplication on functions from the space
$L^{\infty }({\tilde \Gamma }_X,P_m,{\bf C})$, where $\bar H:=L^2({\tilde 
\Gamma }_X,P_m,{\bf C})$.
If measures $P_m$ and $P_{m'}$ are singular, then 
\par $(IV.ii)$ either 
$\sup_{(\gamma \in {\tilde \Gamma }_X)}|P_{m'}(d\gamma )/P_m(d\gamma )|=\infty $
or $\sup_{(\gamma \in {\tilde \Gamma }_X)}|P_m(d\gamma )/P_{m'}
(d\gamma )|=\infty $, where $P_{m'}(d\gamma )/P_m(d\gamma ):=
\lim_{(P_m(B)\to 0, \infty >P_m(B)>0, \gamma \in B)}P_{m'}(B)/P_m(B) \in
[0,\infty ]$, $[0,\infty ]:=([0,\infty )\cup \{ \infty \} )$,
$[0,\infty ):= \{ x: x\in {\bf R}, 0\le x \} $, $B\in Bf({\tilde \Gamma }_X)$.
In view of the existence of the intertwining operator $T$ of $U_m$
with $U_{m'}$ there exists an isomorphism of Hilbert spaces
$\tau : L^2({\tilde \Gamma }_X,P_m,{\bf
C})\to L^2({\tilde \Gamma }_X,P_{m'},{\bf C})$, which has a continuous extension
to an isomorphism of Banach spaces
$\tau : L^{\infty }({\tilde \Gamma }_X,P_m,{\bf C})
\to L^{\infty }({\tilde \Gamma }_X,P_{m'},{\bf C})$ due to condition
$(IV.i)$. On the other hand, in view of condition $(IV.ii)$ there
exists a sequence $f_n\in L^2({\tilde \Gamma }_X,P_m,{\bf C})$
such that $C_1a_n\le b_n\le C_2a_n$ for each $n\in \bf N$
and $\lim_{n\to \infty }c_n<\infty $ and $\lim_{n\to \infty }d_n=\infty $,
where $C_1$ and $C_2$ are positive constants,
$a_n:=\| f_n\|_{L^2({\tilde \Gamma }_X,P_m,{\bf C})}$,
$b_n:=\| \tau f_n\|_{L^2({\tilde \Gamma }_X,P_{m'},{\bf C})}$,
$c_n:=\| f_n\|_{L^{\infty }({\tilde \Gamma }_X,P_m,{\bf C})}$,
$d_n:=\| \tau f_n\|_{L^{\infty }({\tilde \Gamma }_X,P_{m'},{\bf C})}$,
since there are sequences $\{ y_n: 0<y_n<\infty , n\in {\bf N} \} $
such that $\sum_ny_n^{-2}<\infty $, but $\sum_ny_n^{-1}=\infty $. This means
that singularity of $P_m$ with $P_{m'}$ leads to the contradiction,
consequently, $P_m$ and $P_{m'}$ are equivalent.
The cases $(\beta )$ are proved analogously with $\mu $ instead of $P_m$ and 
$G$ instead of ${\tilde \Gamma }_X$ due to theorem 3.8.
\par {\bf 3.12. Corollary.} {\it $(\alpha )$. If
$U^q_m$ and $U^{q'}_{m'}$ are 
equivalent as unitary representations, then $P_m$ and $P_{m'}$ are 
equivalent as measures. 
\par $(\beta )$. If $T^{\mu }$ and $T^{\mu '} $ are equivalent
as unitary representations, then $\mu $ and $\mu '$ are equivalent as measures.}
\par {\bf 3.13. Theorem.} {\it $(\alpha )$.
If $P_m$ and $P_{m'}$ are equivalent,
$n=n'$, unitary representations $q$ and $q'$ of $\Sigma _n$ 
and $\Sigma _{n'}$ are
equivalent (in the case $Diff^t(M)$ of the real manifold $M$ 
with the additional condition $dim_{\bf R}M>1$; $H^q=\{ 0\} $
and $q=I$
for $X=G$ or for infinite-dimensional manifold $X=M$ over the field 
$\bf L$). Then the unitary
representations $U^q_m$ and $U^{q'}_{m'}$ are equivalent.
\par $(\beta )$. If $\mu $ and $\mu '$ from theorem 3.8 are equivalent
quasi-invaraint measures on $G$ relative to $G'$, then the regular unitary
representations $T^{\mu }$ and $T^{\mu '}$ are equivalent.}
\par {\bf Proof.} The cases $(\alpha )$ for $X=M$ infinite-dimensional
over the field $\bf L$ 
or $X=G$ and $(\beta )$ follow from the fact that $\tau : L^2(Z,\mu ,{\bf C})
\to L^2(Z, \mu ', {\bf C})$ given by the following formula
$(\tau f)(x)=(\mu (dx)/\mu '(dx))^{1/2}f(x)$ is the linear topological
isomorphism and the intertwining operator of two regular representations
in these Hilbert spaces, where either $Z={\tilde \Gamma }_X$ with
$\mu =P_m$ or $Z=G$ with a quasi-invariant measure $\mu $ relative to $G'$
respectively.
\par It remains only to consider the case of the non-Archimedean manifold
$X=M$ with $dim_{\bf L}M<\infty $, since the case of real $M$
was proved in \S 4 \cite{vgg}. The measure $m$ on $M$ is supposed to be the
restriction of the Haar measure from $\bf L^n$ on $M$ (see \S 2.9).
Let $Diff^t(X,m)$ be a subgroup of $G'=Diff^t(X)$ consisting of diffeomorphisms
$\psi $, for which $\rho _m(\psi ,x)=1$ for each $x\in M$, where $1\le t\le
\infty $. In the case
of $X=M=\bf R$ we have $Diff^t(X,m)=\{ e\} $, but in the non-Archimedean case
each $\psi \in G'$ with $sup_{x\in M}
|\psi '(x)-I|<1$ belongs to $Diff^t(X,m)$. 
For example, if a countable family of disjoint balls $B_j:=B({\bf L^n},x_j,r)$
with $j\in \Upsilon \subset \bf N$ of radius $0<r<\infty $
is contained in $M$ and if $\psi \in G'$ is such that $\psi (B_j)=
B_{\zeta (j)}$  for each $j\in \Upsilon $, $\psi |_{B_j}(x_j+z)=
x_{\zeta (j)}+\phi _j(z)$ for each $z\in B({\bf L^n},0,r)$, $\psi |_{(M\setminus
\bigcup_jB_j)}=id$, $\phi _j: B({\bf L^n},0,r)\to B({\bf L^n},0,r)$
are diffeomorphisms with $\sup_{(x\in B({\bf L^n},0,r))} |\phi '_j(x)-I|
<1$ for each $j\in \bf N$, where $\zeta : \Upsilon \to \Upsilon $
is a bijection, then $\psi \in Diff^t(M,m)$,
since $B({\bf L^n},0,r)$ are clopen in $\bf L^n$
and the valuation group $\{ |x|_{\bf L}: 0\ne x\in {\bf L} \} $
is discrete in $(0,\infty )$.
Therefore, in the non-Archimedean 
case there is not any restriction on $dim_{\bf L}M$ from below. If $A\subset
\Gamma _X$ is invariant by $mod (P_m )$ subset of $\Gamma _X$ and $P_m(A)>0$,
then $P_m(A)=1$. Indeed, if $P_m(\Gamma _X\setminus A)>0$, then
there exists $\psi \in Diff^t(X,m)$ such that $P_m((\Gamma _X\setminus A)\cap
\psi (A))>0$, since $m$ is quasi-invariant relative to $G'$ with the continuous
quasi-invariant factor $\rho _m(h,x)$ by $(h,x) \in G'\times \Gamma _X$
and $P_m((hA)\cap B)$ is the continuous function by $h\in G'$ for each
$A$ and $B\in Af(\Gamma _X,P_m)$, where $Af(\Gamma _X,P_m)$ denotes the
completion of the Borel $\sigma $-field $Bf(\Gamma _X)$ by the ergodic
measure $P_m$.
In view of the invariance of $A$ we have $P_m((\Gamma _X\setminus A)\cap A)>0$,
which is a contradiction, hence $P_m(\Gamma _X\setminus A)=0$.
\par The restriction of the regular unitary representation
$U_m|_{Diff^t(X,m)}$ is given by the following formula:
\par $(U_m(\psi )f)(\gamma )=f(\psi ^{-1}\gamma )$ for each $\psi \in
Diff^t(X,m)$. Then $f_0(\gamma )=1$ for each $\gamma \in \Gamma _X$
is the unique vector in $L^2(\Gamma _X,P_m,{\bf C})$ such that
${\bf C}f_0$ is invariant relative to $U_m(\psi )$ for each $\psi
\in Diff^t(M,m)$. The Poisson measure $P_m$ can be considered with a 
parameter $\lambda >0$, that is with $\lambda m$ instead of $m$.
Let $u_m(\psi )$ be a spherical function given by the following formula:
$u_m(\psi )=(U_m(\psi )f_0,f_0)$, where $(*,*)$ denotes the scalar
product in $L^2(\Gamma _X,P_m,{\bf C})$. Then $u_m(\psi )=exp(\int_X
(\rho ^{1/2}_m(\psi ,x)-1)m(dx))$, since $u_m(\psi )=
\int_{\Gamma _X}(\prod_{x\in \gamma }\rho ^{1/2}(\psi ,x))P_m(d\gamma )$
and for $supp(\psi )\subset Y$ with $m(Y)<\infty $ we have
$u_m(\psi )=\sum_{n=0}^{\infty }\int_{B^n_Y}(\prod_{x\in \gamma }
\rho ^{1/2}(\psi ,x))P_m|_{B^n_Y}(d\gamma )$. Therefore, we get the following
theorem.
\par {\bf 3.14. Theorem.} {\it $(\alpha )$.
The representations $U_{\lambda _1m}$
and $U_{\lambda _2m}$ of $Diff^t(M)$ (with the restriction 
$dim_{\bf R}M>1$
for $M$ over $\bf R$ and $0<dim_{\bf L}M$ in the 
non-Archimedean case) for $\lambda _1\ne \lambda _2$
are inequivalent.
\par $(\beta )$. The representations $U_{\lambda _1m}$ and $U_{\lambda _2m}$
of $G'$ in $L^2({\tilde \Gamma }_G,P_{\lambda _jm},{\bf C})$ with $j=1,2$
respectively are inequivalent for $\lambda _1\ne \lambda _2$.}
\par {\bf Proof.} $(\alpha )$. For $dim_{\bf L}M<\infty $ 
this follows from the fact
$u_{\lambda _1m}\ne u_{\lambda _2m}$.
\par $(\beta )$. In view of the Kakutani theorem \cite{dal} 
two Poisson measures $P_{\lambda _1m}$
and $P_{\lambda _2m}$ are singular, hence by theorem 3.10 representations
$U_{\lambda _1m}$ and $U_{\lambda _2m}$ are inequivalent.
\par {\bf 3.15. Note.} By the given above representations it is possible to
produce new with the help of the following construction.
Let $G'$ be a group acting from the left on a $C^{\infty }$-manifold $Y$
(or analytic with disjoint charts in the non-Archimedean case over the 
local field $\bf L$)
such that on $Y$ is given a $\sigma $-additive $\sigma $-finite
quasi-invariant non-negative measure $\mu $ with a continuous
quasi-invariant factor $\rho _{\mu }(\psi ,y)$ by $(\psi ,y)\in G'\times Y$.
In the real case let us consider a space 
${\sf F}(Y)$ of generalised functions on $Y$. For example, if there is a 
unitary regular representation $T$ of $G'$ in $L^2(Y,\mu ,{\bf C})$, then
${\sf F}(Y)$ is a space of continuous linear functionals on
$L^2(Y,\mu ,{\bf C})$, hence ${\sf F}(Y)$ is isomorphic with
$L^2(Y,\mu ,{\bf C})$. Let $\nu $ be a measure on ${\sf F}(Y)$
given with the help of its characteristic function
$\int_{{\sf F}(Y)}exp(i<F,f>)\nu (dF)=exp(-\| f\| ^2/2)$,
where $\| *\| $ is a norm in $L^2(Y,\mu ,{\bf C})$. Such $\nu $ is called
the standard Gaussian measure in ${\sf F}(Y)$. Then a new representation
${\tilde U}:=EXP_{\beta }T$ is given by the following formula:
\par $({\tilde U}(\psi )\Phi )(F):=exp(i<F,\beta (\psi )>)\Phi (T^*(\psi )F)$,
where $<T^*(\psi )F,f>=<F,T(\psi )f>$, $<F,f>$ is a value of a functional
$F\in {\sf F}(Y)$ on a function $f\in L^2(Y,\mu ,{\bf C})$, $\beta $ is a
1-cocycle such that $[\beta (\psi )](y):=\rho ^{1/2}_{\mu }(\psi ,y)-1$,
$T(\psi ):=\rho ^{1/2}(\psi ,y)f(\psi ^{-1}y)$ for each $\psi \in G'$,
$\Phi \in L^2({\sf F},\nu ,{\bf C})$, $T^*(\psi )=T^{-1}(\psi )=
T(\psi ^{-1})$. If to substitute $\beta $ on $s\beta $,
where $s\in \bf R$, then it produces the one-parameter family ${\tilde U}_s:=
EXP_{s\beta }T$. There is an equality $\lim_{\psi \to e}\beta (\psi ,y)=0$
for each $y\in Y$. When $Y=G$ or $Y={\tilde \Gamma }_X$ and $\mu $
is as in theorems 3.8 or 3.9, then for $s\ne 0$ representation ${\tilde U}_s$
is irreducible as follows from the proof of theorems 3.8, 3.9, since
the linear span of non-linear functionals $exp(i<F,\beta (\psi )>)$
is dense in $L^2({\sf F}(Y),\nu ,{\bf C})$ (see also \S 4 in \cite{vgg}).
For $X=M$ with $dim_{\bf L}M<\infty $ and $s\ne 0$
the representations ${\tilde U}_s$
and $U_{s^2\mu }$ in $L^2(\Gamma _M,P_{\mu },{\bf C})$
are equivalent, since $({\tilde U}_s{\psi }\Phi _0,\Phi _0)=exp(\int_M
(\rho ^{1/2}(\psi ,x)-1)\mu (dx))=u_{s^2}(\psi )$, where $\Phi _0(F)=1$
for each $F\in {\sf F}(\Gamma _M)$, $\Phi _0\in L^2({\sf F}(\Gamma _M),
\nu ,{\bf C})$.
\par {\bf 3.16. Note.}It follows from \cite{dal,lupr210,sko}, that
on $Z$ there are infinite families of orthogonal measures, restrictions of 
which on $U$ are
quasi-invariant relative to $G'$ and have continuous quasi-invariance
factor on $G'\times X$, where either $X=M$ or $X=G$ respectively,
$W'U\subset V$, $W'$ is open in $G'$ and $U$ is open in $Z$. 
Due to the general procedures of construction of measures on 
$X$ outlined above
on infinite-dimensional $M$ over the corresponding field
or on $G$ there are infinite families of orthogonal and as well
singular measures, since measures on these infinite-dimensional manifolds
$M$ and $G$ are induced from the corresponding Banach spaces $Z$ due to the
local diffeomorphisms $A: W\to V$, where $W$ is open in $M$ or $G$
and $V$ is open in $Z$. 
Therefore, the last two theorems show that there exists an infinite
family of non-equivalent unitary representations of $G'$ for $X=G$
and also for $G'=Diff^t(M)$ for the infinite-dimensional manifold $M$
over the corresponding field, since in these cases on $G$ and $M$ 
there exist infinite families of orhtogonal measures.
The unitary group $U(l_2)$ of the standard Hilbert space $l_2$ over
$\bf C$ has the topological density $\sf c$, when $U(l_2)$ is in its standard
topolgy induced by the operator norm in the space of linear bounded operators
$L(l_2)$ on $l_2$, since $l_2^{\bf N}$ in the box topology has a density
$\aleph _0^{\aleph _0}=\sf c$. When $U(l_2)$ is considered as a topological
space in its strong topology \cite{fell}, then its topological density is
$\aleph _0$, since $l_2^{\bf N}$ in the product Tychonoff
topology has density $\aleph _0$ \cite{eng}. Therefore,
the cardinality of distinct unitary strongly continuous representations
$T: G'\to U(l_2)$ for topological group with density $\aleph _0$
do not exceed $\sf c$, since ${\sf c}^{\aleph _0}=\sf c$ \cite{eng}.
This is important difference of the theory of such non-locally compact
topological groups with the theory of compact groups. In the latter
case all irreducible unitary representations arise as irreducible
components of the regular representation associated with the Haar measure, but
for the considered here cases of groups this is not true, since there are
infinite families of non-equivalent unitary representations on such groups.
There are considered $M$ and $G$ and ${\tilde \Gamma }_X$ with countable bases
of topology and real-valued measures. The family $\Psi $
of distinct $\sigma $-additive Borel measures on these spaces 
have the cardinality $card(\Psi )=card({\bf R}^{\bf N})=card({\bf R})=:\sf c$.
In view of theorems 3.10, 3.13, 3.14 and  
the criteria of orthogonality and singularity of measures
on infinite-dimensional spaces (using weak distributions, product measures,
Kakutani theorem and its non-Archimedean analog \cite{dal,lupr210}
and the construction of measures on $G$ or ${\tilde \Gamma }_X$ with the
help of local diffeomorphisms of open subsets in these spaces and neigbourhoods
of zero in the corresponding Banach spaces as in \S 2.9) 
there exist families $\Psi _s$
of singular and $\Psi _o$ of orthogonal measures such that ${\sf c}\le
card(\Psi _o)\le card(\Psi _s)\le card(\Psi )=\sf c$, hence there are
$\sf c$ non-equivalent unitary representations $U_m$ of $G'$ 
in $L^2({\tilde \Gamma }_X,P_m,{\bf C})$ and also $\sf c$ inequivalent unitary
representations of $G'$ in $L^2(G,\mu ,{\bf C})$, which were considered above,
since $card(\Psi _s)=\sf c$.
\par {\bf 3.17. Theorem.} {\it There exist 
Abelian non-locally compact Banach-Lie 
groups $G$ with quasi-invariant measures $\mu $ on $G$
relative to dense subgroups $G'$ such that the associated regular unitary
representations $T^{\mu }$ of $G'$ are irreducible, 
each one-parameter subgroup $S$ of $G$ is compact and a projection
of $\mu $ on each one-parameter subgroup $S$
is equivalent with the Haar measure on $S$.}
\par {\bf Proof.} Let $l_{2,b}$ be a Hilbert space over $\bf R$
of elements $x=(x_j: j\in {\bf N}, x_j\in {\bf R} )$ such that
$\| x\| ^2_{l_{2,b}}:=\sum_j|x_jj^b|^2<\infty $. In particular
$l_2=l_{2,0}$. These spaces have standard orthonormal bases
$e_k:=(0,...,0,1,0,... )$ with $1$ on the $k$-th place, $k\in \bf N$.
For a local field $\bf L$ let $c_{0,b}$ be a Banach space of elements
$x=(x_j: j\in {\bf N}, x_j\in {\bf L} )$ such that $\| x\|:=
\max_j|x_j|_{\bf L}p^{jb}<\infty $ and $\lim_{j\to \infty }|x_j|p^{jb}=0$,
where $p$ is a prime number such that $\bf L$ is a finite algebraic extension 
of the field of $p$-adic numbers $\bf Q_p$. In particular $c_{0,0}=c_0$.
If $b>1$ (with $p^{-b}\in {\Gamma '}_{\bf L}$ in the non-Archimedean case)
then the embeddings $J_b: l_{2,b}\hookrightarrow l_2$
and $S_b: c_{0,b}\hookrightarrow c_0$ are of trace class:
$J_be_k=a_ke_k$ with $a_k=k^{-b}\in \bf R$ and $\sum_k|a_k|<\infty $,
$S_be_k=v_ke_k$ with $v_k\in \bf L$ and $\sum_k|v_k|<\infty $,
where $|v_k|=p^{-jb}$. On $l_2$ and $c_0$ there exist a Gaussian measure
$\lambda $ and a non-Archimedean analog $\eta $
of a Gaussian measure quasi-invariant relative to
$l_{2,b}$ and $c_{0,b}$ respectively such that their projections $\lambda _k$
and $\eta _k$ on one-dimensional subspaces ${\bf R}e_k$ and ${\bf L}e_k$
are the following: $\lambda _k(dx_k)=C_kexp(-x_k^2s_k^2)w(dx_k)$
and $\eta _k(dy_k)=F_kexp(-|y_k|^2p^{2k})v(dy_k)$, 
where $w$ and $v$ are the Lebesgue
and the Haar measures on $\bf R$ and $\bf L$ respectively
such that $w([0,1])=1$, $v(B({\bf L},0,1))=1$, $s_k=k^{b'}$ for each $k\in \bf
N$ with $1<b'<b$, $C_k>0$, $F_k>0$, $\lambda _k({\bf R})=1$, 
$\eta _k({\bf L})=1$ (see \cite{dal,lupr210}).
Consider an additive discrete subgroup $E$ of $l_{2,b}$ consisting of elements
$x\in l_{2,b}$ such that $x_j=n_je_j$ for each $j\in \bf N$, where 
$n_j\in \bf Z$. Then $l_{2,b}/E=:H_b$ and $l_2/E=:H $ are the additive groups.
The measures $\lambda $ and $\eta $ induce measures $\mu $ on $H$ and 
$\nu $ on $B(c_0,0,1):=\{ x\in c_0: \| x\| \le 1 \}=:B$. Then $\mu $
is quasi-invariant relative to $H_b$ and $\nu $ is quasi-invariant relative to
$B(c_{0,b},0,1)$ with continuous quasi-invariance factors
such that $\mu (H)=1$ and $\nu (B)>0$. 
\par Let $L_n:=sp_{\bf R}(e_1,...,e_n)$ and $E_n:=L_n\cap E$, so the latter
is a discrete subgroup of $L_n$ and $L_n/E_n=:V_n$ is a closed subgroup of
$H_b$. Hence a projection $\pi _n: l_{2,b}\to L_n$, which has a continuous
extension $\pi _n: l_2\to L_n$ induces a quotient mapping ${\bar \pi }_n:
H_b\to V_n$ with a continuous extension ${\bar \pi }_n: H\to V_n$ for each
$n\in \bf N$. Therefore, the measure $\mu $ on $H$ induces a measure
$\mu _n$ on $V_n$ such that $\mu _n(A):=\mu ({\bar \pi }_n^{-1}(A))$
for each $A\in Bf(V_n)$. In view of the equality $\lim_{n\to \infty }
\rho _{\mu _n}({\bar \pi }_n(\psi ), {\bar \pi }_n(x))=\rho _{\mu }(\psi ,x)$
for each $\psi \in H_b$ and $x\in H$ it follows that $\rho _{\mu }(\psi ,x)=
\lim_{n\to \infty }(\sum_{z\in E_n}exp\{ \sum_{l=1}^n[2(\psi +z,e_l)(x,e_l)-
(\psi +z,e_l)^2]s_l^2 \} )$ $(\sum_{z\in E_n}exp\{
\sum_{l=1}^n[2(z,e_l)(x,e_l)-(z,e_l)^2]s_l^2 \} )^{-1}$, since 
$(\pi _n(x),e_l)=(x,e_l)$ for each $x\in L_n$ with $n\ge l$.
The Hilbert space $L^2(H,\mu ,{\bf C})$ is isomorphic with
a subspace $\{ f: f\in L^2(l_2,\lambda ,{\bf C}); f(x+z)=f(x)$
$\lambda $-a.e. for each $z\in E \} $. Since $sp_{\bf C}
\{ \rho ^{1/2}_{\lambda }(\psi ,x)=:\phi (x)| \psi \in l_{2,b} \} $ is dense
in $L^2(l_2,\lambda ,{\bf C})$, then $sp_{\bf C} \{ \rho ^{1/2}_{\mu }
=:\phi (x)| \psi \in H_b \} $ is dense in $L^2(H,\mu ,{\bf C})$.
Repeating the proof of theorem
3.8 for these groups we get that their regular unitary representations
are irreducible. That is representations $T$ of $l_{2,b}$,
$H_b$, $c_{0,b}$ and $B(c_{0,b},0,1)$ in the corresponding spaces 
$L^2(l_2,\lambda ,{\bf C})$, $L^2(H,\mu ,{\bf C})$,
$L^2(c_0,\eta ,{\bf C})$ and $L^2(B,\nu ,{\bf C})$ (the first case
was also considered more generally for additive groups of locally convex
spaces in \cite{banas,gelvil}).
These groups are Banach-Lie and Abelian, moreover, each one-parameter
subgroup of $H_b$ and of $B(c_{0,b},0,1)$ over $\bf R$ and $\bf L$
respectively is compact. The projections of $\mu $ and $\nu $
on one-parameter subgroups are equivalent with the Haar measures 
on them. Certainly, $H_b$ and $B(c_{0,b},0,1)$ are not locally compact,
since $T_eH_b$ and $T_eB(c_{0,b},0,1)$ are infinite-dimensional
Banach spaces over $\bf R$ and $\bf L$ respectively.
\par {\bf 3.18. Note.} Regular representations $U_m$
of the groups $l_{2,b}$, $H_b$, $c_{0,b}$ or $B(c_{0,b},0,1)$
from the proof of theorem 3.17 in the space 
$L^2({\tilde \Gamma }_X,P_m,{\bf C})$ with $X=l_2$, $H$, $c_0$ or $B$ and
$m=\lambda $, $\mu $, $\eta $ or $\nu $ respectively are reducible, since
$D^nL_h=0$ for each $n>1$ and $f_0$ is not cyclic for $U_m$
(see the proof of theorem 3.9).


\begin{thebibliography}{99}
\bibitem{banas} W. Banaszczyk. Additive subgroups of topological 
vector spaces (Berlin: Springer, 1991).
\bibitem{bao} D. Bao, J. Lafontaine, T. Ratiu. On a non-linear 
equation related to the geometry of the difeomorphism group.
Pacif. J. Mathem. {\bf 158} (1993), 223-242.
\bibitem{boug} N. Bourbaki. Lie groups and algebras (Moscow: Mir, 1976).
\bibitem{boui} N. Bourbaki. Integration. Chapters 1-9 (Moscow: Nauka, 1970
and 1977).
\bibitem{dal} Yu.L. Dalecky, S.V. Fomin. Measures and differential
equations in infinite-dimensional space (Kluwer: Dordrecht, The Netherlands,
1991).
\bibitem{ebi} D.G. Ebin, J. Marsden. Groups of diffeomorphisms and the motion of
incompressible fluid. Ann. Math. {\bf 92} (1970), 102-163.
\bibitem{eng} R. Engelking. General topology (Moscow: Mir,1986).
\bibitem{fed} H.Federer. Geometric measure theory(Berlin:Springer-Verlag,
    1969).
\bibitem{fell} J.M.G. Fell, R.S. Doran. Representations of $*$-algebras,
locally compact groups, and Banach $*$-algebraic bundles (Acad. Pr.:
Boston, 1988).
\bibitem{gelvil} I.M. Gelfand, N.Ya. Vilenkin. Generalized
functions. v. 4. Some applications of harmonic analysis. Rigged 
Hilbert spaces. (Moscow: Nauka, 1961).
\bibitem{hew} E. Hewitt and K.A. Ross. Abstract harmonic analysis. Second
edition (Berlin: Springer-Verlag, 1979).
\bibitem{hira} T. Hirai. Irreducible unitary representations of the group of
diffeomorphisms of a non-compact manifold. J. Math. Kyoto Univ.
{\bf 33} (1993), 827-864.
\bibitem{kling} W. Klingenberg. Riemannian geometry (Walter de Gruyter:
Berlin, 1982).
\bibitem{kos} A.V. Kosyak. Irreducible Gaussian representations
of the group of the interval and circle diffeomorphisms.
J. Funct. Anal. {\bf 125}(1994), 493-547.
\bibitem{kuo} H.-H. Kuo. Gaussian measures in Banach spaces (Springer,
Berlin, 1975).
\bibitem{luummr} S.V. Ludkovsky. Measurability of representations
of infinite-dimensional groups. {\bf 51} (1996), 205-206 (N 3).
\bibitem{luumn96} S.V. Ludkovsky. Measures on groups of diffeomorphisms of
non-Archimedean Banach manifolds, Usp. Mat. Nauk.
{\bf 51}(1996), 169-170 (N 2).
\bibitem{lutmf} S.V. Ludkovsky. Measures on groups of diffeomorphisms of
non-Archimedean manifolds, representations of groups and their applications.
Theoret. i Mathem. Phys., 1999.
\bibitem{luumnad} S.V. Ludkovsky. Quasi-invariant measures on 
non-Archimedean semigroups of loops.
Usp. Mat. Nauk, {\bf 53} (1998), 203-204 (N 3).
\bibitem{luseamb} S.V. Ludkovsky. Irreducible unitary representations of 
non-Archimedean groups of diffeomorphisms.
Southeast Asian Bulletin of Mathematics (Hong Kong).
{\bf 22} (1998), 419-436.
\bibitem{lurimut} S.V. Ludkovsky.
Irreducible unitary representations
of a diffeomorphisms group of an infinite-dimensional real manifold.
 Rendiconti dell'Istituto di Matematica dell'Universit\`a
di Trieste. Nuova Serie, 
(in English)  {\bf 29} (1998), 22 pages.
\bibitem{luum985} S.V. Ludkovsky. Embedding of a non-Archimedean 
Banach manifold into the corresponding Banach space. Usp. Mat. Nauk.,
{\bf 53} (1998), N 5. 
\bibitem{luihlgr} S.V. Ludkovsky. Gaussian quasi-invariant measures
on loop groups and semigroups of real manifolds and their representations.
IHES, Bures-sur-Yvette, France, preprint {\bf IHES/M/97/95}.
\bibitem{luihlgna} S.V. Ludkovsky. Quasi-invariant measures
on non-Archimedean groups and semigroups of loops and paths,
their representations. {\bf IHES/M/98/36}.
\bibitem{luihgdsw} S.V. Ludkovsky. Quasi-invariant measures
on groups of diffeomorphisms of Schwarz class of smoothness
for real manifolds. IHES/M/97/96.
\bibitem{lupr180} S.V. Ludkovsky. Representations and structure of groups of
diffeomorphisms of non-Archimedean Banach manifolds. I, II. 
Intern. Centre for Theoret. Phys. 
Trieste, Italy. Preprints
N {\bf IC/96/180,181}, September 1996 (http://www.ictp.trieste.it).
\bibitem{lupr202} S.V. Ludkovsky. Quasi-invariant measures on a group
of diffeomorphisms of an infinite-dimensional Hilbert manifold and
its representations.
ICTP. {\bf IC/96/202}, October 1996.
\bibitem{lupr210} S.V. Ludkovsky. Quasi-invariant and pseudo-differentiable
measures on a non-Archimedean Banach space.
ICTP. {\bf IC/96/210}, October 1996.
\bibitem{lupr215} S.V. Ludkovsky. Quasi-invariant
measures on a non-Archimedean group
of diffeomorphisms and on a Banach manifold. ICTP.
{\bf IC/96/215}, October, 1996.
\bibitem{lupr218} S.V. Ludkovsky. Quasi-invariant
measures on groups of diffeomorphisms of real Banach 
manifolds. ICTP. {\bf IC/96/218}, October, 1996.
\bibitem{nai} M.A. Naimark. Normed rings (Moscow: Nauka, 1968).
\bibitem{nari} L. Narici, E. Beckenstein. Topological vector spaces.
(Marcel Dekker Inc.: New York, 1985).
\bibitem{nere} Yu.A. Neretin. Representations of the Virasoro algebra and
    affine algebras. in: Itogi Nauki i Tech. Ser. Sovr. Probl.
    Math. Fund. Napravl(Moscow: Nauka) {\bf 22}(1988), 163-230.
\bibitem{roo} A.C.M. van Rooij. Non-Archimedean functional analysis
(Marcel Dekker Inc.: New York, 1978).
\bibitem{shav} E.T. Shavgulidze. About one measure quasi-invariant relative
    to an action of a diffeomorphisms group of a finite-dimensional
    manifold. Dokl. Akad. Nauk SSSR. {\bf 303}(1988), 811-814.
\bibitem{shim} H. Shimomura. Poisson measures on the configuration space
and unitary representations of the group of diffeomorphisms.
J. Math. Kyoto Univ. {\bf 34} (1994), 599-614.
\bibitem{sko} A.V. Skorohod. Integration in the Hilbert space (Moscow:
    Nauka, 1975).
\bibitem{vgg} A.M. Vershik, I.M. Gelfand, M.I. Graev.
Representations of the group of diffeomorphisms.
Usp. Mat. Nauk. {\bf 30} (1975), 3-50.
\end{thebibliography}
\end{document}